\documentclass[11pt,a4paper,reqno]{amsart}

\usepackage{epsfig}
\usepackage{amsfonts}
\usepackage{amsmath}
\usepackage{amssymb}
\usepackage{amsthm}
\usepackage{color}
\usepackage{enumerate}
\usepackage[T1]{fontenc}
\usepackage[latin1]{inputenc}
\usepackage{ae,aecompl}
\usepackage{graphicx}
\usepackage{nicefrac}
\usepackage[includehead,includefoot,margin=25mm]{geometry}

\newtheorem{theorem}{Theorem}[section]
\newtheorem{corollary}[theorem]{Corollary}
\newtheorem{lemma}[theorem]{Lemma}
\newtheorem{proposition}[theorem]{Proposition}

{\theoremstyle{remark}
  \newtheorem{remark}[theorem]{Remark}}
{\theoremstyle{definition}
  \newtheorem{definition}[theorem]{Definition}
  
  \newtheorem{notation}[theorem]{Notation}
  \newtheorem{example}[theorem]{Example}
  
  \newtheorem{sit}[theorem]{}
}

\newcommand{\PP}[0]{\ensuremath{\mathbb{P}}}
\newcommand{\CC}[0]{\ensuremath{\mathbb{C}}}
\newcommand{\ZZ}[0]{\ensuremath{\mathbb{Z}}}
\newcommand{\AF}[0]{\ensuremath{\mathbb{A}}}

\newcommand{\QQ}[0]{\ensuremath{\mathbb{Q}}}
\newcommand{\TT}[0]{\ensuremath{\mathbb{T}}}
\newcommand{\KK}[0]{\ensuremath{\mathbf{k}}}
\newcommand{\OO}[0]{\ensuremath{\mathcal{O}}}
\newcommand{\DD}[0]{\ensuremath{\mathfrak{D}}}
\newcommand{\EE}[0]{\ensuremath{\mathcal{E}}}

\newcommand{\fract}[0]{\ensuremath{\mathrm{Frac}}}
\newcommand{\LND}[0]{\ensuremath{\mathrm{LND}}}
\newcommand{\nil}[0]{\ensuremath{\mathrm{Nil}}}

\newcommand{\spec}[0]{\ensuremath{\mathrm{Spec}}}

\newcommand{\ML}[0]{\ensuremath{\mathrm{ML}}}

\newcommand{\divi}[0]{\ensuremath{\mathrm{div}}}

\newcommand{\ord}[0]{\ensuremath{\mathrm{ord}}}
\newcommand{\pol}[0]{\ensuremath{\mathrm{Pol}}}
\newcommand{\cpl}[0]{\ensuremath{\mathrm{CPL}}}
\newcommand{\rank}[0]{\ensuremath{\mathrm{rank}}}
\newcommand{\cone}[0]{\ensuremath{\mathrm{cone}}}
\newcommand{\homo}[0]{\ensuremath{\mathrm{Hom}}}

\begin{document}

\title[Affine $\TT$-varieties and locally nilpotent
derivations]{Affine $\TT$-varieties of complexity one and locally
  nilpotent derivations}

\author{Alvaro Liendo}

\address{Universit\'e Grenoble I, Institut Fourier, UMR 5582 CNRS-UJF,
  BP 74, 38402 St.\ Martin d'H\`eres c\'edex, France}

\email{alvaro.liendo@ujf-grenoble.fr}

\date{\today}

\thanks{ \mbox{\hspace{11pt}}{\it 2000 Mathematics Subject
    Classification}:
  14R05, 14R20, 13N15, 14M25.\\
  \mbox{\hspace{11pt}}{\it Key words}: torus action, $\KK_+$-action,
  locally nilpotent derivations, affine varieties, Makar-Limanov
  invariant}

\begin{abstract}
  Let $X=\spec\, A$ be a normal affine variety over an algebraically
  closed field $\KK$ of characteristic 0 endowed with an effective
  action of a torus $\TT$ of dimension $n$. Let also $\partial$ be a
  homogeneous locally nilpotent derivation on the normal affine
  $\ZZ^n$-graded domain $A$, so that $\partial$ generates a
  $\KK_+$-action on $X$ that is normalized by the $\TT$-action.

  We provide a complete classification of pairs $(X,\partial)$ in two
  cases: for toric varieties ($n=\dim X$) and in the case where
  $n=\dim X-1$. This generalizes previously known results for surfaces
  due to Flenner and Zaidenberg. As an application we compute the
  homogeneous Makar-Limanov invariant of such varieties. In particular
  we exhibit a family of non-rational varieties with trivial
  Makar-Limanov invariant.
\end{abstract}

\maketitle

\setcounter{tocdepth}{3}
\tableofcontents

\section*{Introduction}

Let $\KK$ be an algebraically closed field of characteristic 0. For an
algebraic torus $\TT\simeq(\KK^*)^n$ acting on an algebraic variety
$X$, the complexity of this action is the codimension of the general
orbit. Without loss of generality, we restrict to effective
$\TT$-actions, so the complexity is $\dim X-\dim \TT$. In particular,
a $\TT$-variety of complexity 0 has an open orbit and thus is a toric
variety. It is well known that a $\TT$-action on $X=\spec\,A$ gives
rise to an $M$-grading on $A$, where $M$ is a lattice of rank $n$.

More generally, let $A=\bigoplus_{m\in M}\tilde{A}_m$ be a finitely
generated effectively $M$-graded domain and $K=\fract\,A$. For any
$m\in M$ we let
$$K_m=\left\{f/g\in K\mid f\in \tilde{A}_{m+e}, g\in \tilde{A}_e\right\}\,.$$
Then $\tilde{A}_m\subseteq K_m$, and $\KK \subseteq K_0 \subseteq K$
are field extensions. Letting $\{\mu_1,\ldots,\mu_n\}$ be a basis of
$M$, we fix for every $i=1,\ldots,n$ an element $\chi^{\mu_i}\in
K_{\mu_i}$. For every $m=\sum_i a_i\mu_i$ we have $K_m=\chi^mK_0$,
where $\chi^m=\prod_i(\chi^{\mu_i})^{a_i}$. Thus, without loss of
generality, we assume in the sequel that
$$ A=\bigoplus_{m\in M} A_m\chi^m\subseteq K_0[M], \quad\mbox{where}\quad A_m\subseteq K_0\,,$$
and $K_0[M]$ denotes the semigroup $K_0$-algebra of $M$. In this
setting, the complexity of the $\TT$-action equals the transcendence
degree of $K_0$ over $\KK$. In particular, for a toric variety $X$,
$K_0=\KK$, and $\chi^m$ is just a character of $\TT$ regarded as a
rational function on $X$.

There are well known combinatorial descriptions of normal
$\TT$-varieties. For toric varieties see e.g., \cite{De}, Chapter 1 in
\cite{KKMS}, and \cite{Od}. For complexity 1 case see Chapters 2 and 4
in \cite{KKMS}, and more generally \cite{Ti1,Ti2}. Finally for
arbitrary complexity see \cite{AlHa,AHS}\footnote{In the case of
  complexity 1, the descriptions in \cite{AlHa} and \cite{Ti2} are
  equivalent and agree with the one in \cite[Chapters 2 and 4]{KKMS},
  see \cite[Section 6]{Ti2} and \cite{Vo}.}.

We let $N=\homo(M,\ZZ)$, $M_{\QQ}=M\otimes\QQ$, and
$N_{\QQ}=N\otimes\QQ$. Any affine toric variety can be described via
the weight cone $\sigma^\vee\subseteq M_{\QQ}$ spanned over
$\QQ_{\geq0}$ by all $m\in M$ such that $A_m\neq\{0\}$ or,
alternatively, via the dual cone $\sigma\subseteq N_{\QQ}$. Similarly,
the description of normal affine $\TT$-varieties of complexity 1 due
to Altmann and Hausen deals with a polyhedral cone $\sigma\subseteq
N_{\QQ}$ (dual to the weight cone $\sigma^\vee\subseteq M_{\QQ}$), a
smooth curve $C$, and a divisor $\DD$ on $C$ whose coefficients are
polyhedra in $N_{\QQ}$ having tail cone $\sigma$. The degree $\deg\DD$
is defined as the Minkowski sum of the coefficients of $\DD$ (see
Subsection \ref{comb-des} for precise definitions).

For affine surfaces with a $\CC^*$-action an alternative
description\footnote{Which is actually equivalent, see Example 3.5 in
  \cite{AlHa}.} was proposed in \cite{FlZa1}. This description was
used in \cite{FlZa2} in order to classify all $\CC_+$-actions on
normal $\CC^*$-surfaces. Generalizing this construction, in the
present paper we use the description in \cite{AlHa} to classify normal
affine $\TT$-varieties of complexity 0 or 1 endowed with a
$\KK_+$-action.

A $\KK_+$-action gives rise to a locally nilpotent derivation (LND) on
$A$. To any LND on $A$ we can associate a homogeneous LND which maps
homogeneous elements into homogeneous elements, see Lemma
\ref{hLDN}. A homogeneous LND $\partial$ on $A=\bigoplus_{m\in M}
A_m\chi^m\subseteq K_0[M]$ can be extended to a derivation on
$K_0[M]$. We say that $\partial$ is of \emph{fiber type} if
$\partial(K_0)=0$ and of \emph{horizontal type} otherwise. In
geometric terms, the fact that the LND $\partial$ is homogeneous means
that the corresponding $\KK_+$-action on $X=\spec\, A$ is normalized by the torus
$\TT$.

In Theorem \ref{toric-clasification} we obtain a classification of
homogeneous LNDs on toric varieties. For $\TT$-varieties of complexity
1, such a classification is given in Theorems \ref{fiber-clas} (for
fiber type) and \ref{resultado} (for horizontal type). These theorems
are the main results of the paper. In \cite{Li1} this classification
of homogeneous LNDs of fiber type is generalized to arbitrary
complexity.

We show as a corollary that the equivalence classes of homogeneous
LNDs on the toric variety defined by the cone $\sigma\subseteq
N_{\QQ}$ are in one to one correspondence with the extremal rays of
$\sigma$ (see Corollary \ref{toric-neq}). This is also true for normal
affine $\TT$-varieties of complexity 1 over an affine curve $C$. Over
a projective curve $C$, these classes are in one to one correspondence
with the extremal rays of $\sigma$ disjoint from the polyhedron
$\deg\DD$ (see Remark \ref{feq2}). The classification of homogeneous
LNDs of horizontal type is more involved, see Corollary \ref{hor-eq}.

The Makar-Limanov invariant \cite{ML} is an important tool which
allows, in particular, to distinguish certain varieties from the
affine space. For an algebra $A$, this invariant is defined as the
intersection of the kernels of all locally nilpotent derivations on
$A$. For graded algebras, we introduce a homogeneous version of the
Makar-Limanov invariant. For $\TT$-varieties of complexity 0 and 1 we
give an explicit expression of the latter invariant. The triviality of
the homogeneous Makar-Limanov invariant implies that of the usual
one.

As an application we exhibit in Subsection \ref{non-rat} a family of
non-rational singular varieties with a trivial Makar-Limanov
invariant. These examples (in a preliminary version of our paper)
attracted the attention of V. L. Popov, who proposed in a recent
preprint \cite{Po} yet another family of affine varieties with these
same properties, this time in addition smooth. It is worthwhile
mentioning that generalizing the methods in Subsection~\ref{non-rat}
we obtained a birational characterization of normal affine varieties
with trivial Makar-Limanov invariant \cite{Li1}.

The content of the paper is as follows. In Section 1 we recall the
combinatorial description of $\TT$-varieties due to Altmann and
Hausen, and also some generalities on locally nilpotent derivations
and $\KK_+$-actions. In Sections \ref{toric} and \ref{com1} we obtain
our principal classification results for toric varieties and for
$\TT$-varieties of complexity 1, respectively. The comparison with
previously known results in the surface case is given in subsection
\ref{sur-cas}. Finally in Section \ref{app} we provide the
applications to the Makar-Limanov invariant.

In the entire paper $\KK$ is an algebraically closed field of
characteristic 0, except in Section \ref{toric}, where $\KK$ is not
necessarily algebraically closed.

The author is grateful to Mikhail Zaidenberg for posing the problem
and permanent encouragement, and to Dimitri Timashev for useful
discussions. We thanks also Vladimir Popov for kindly communicating to
us his preprint \cite{Po}.

\section{Preliminaries}
\label{pre}

\subsection{Combinatorial description of $\TT$-varieties}
\label{comb-des}

Let $N$ be a lattice of rank $n$ and $M=\homo(N,\ZZ)$ be its dual
lattice. We also let $N_{\QQ}=N\otimes\QQ$, $M_{\QQ}=M\otimes\QQ$, and
we consider the natural duality $M_{\QQ}\times N_{\QQ}\rightarrow
\QQ$, $(m,p)\mapsto \langle m,p\rangle$.

Let $\TT=\spec\,\KK[M]$ be the corresponding $n$-dimensional algebraic
torus associated to $M$. Thus $M$ is the character lattice of $\TT$
and $N$ is the lattice of 1-parameter subgroups. It is customary to
write the character associated to a lattice vector $m\in M$ as
$\chi^m$, so that $\chi^m$ is the comorphism of the morphism
$\KK[t]\rightarrow\KK[M],\ t\mapsto m$ \cite{Od}.

Let $X=\spec\,A$ be an affine $\TT$-variety. It is well known that the
morphism $A\rightarrow A\otimes \KK[M]$ induces an $M$-grading on $A$
and, conversely, every $M$-grading on $A$ arises in this
way. Furthermore, a $\TT$-action is effective if an only if the
corresponding $M$-grading is effective\footnote{We say that an
  $M$-graded algebra $A$ is effectively graded by $M$ if the set
  $\{m\in M\mid A_m\neq 0\}$ is not contained in a proper sublattice of
  $M$.}.

Let $A=\bigoplus_{m\in M}A_m\chi^{m}$ be a finitely generated
effectively $M$-graded domain. The \emph{weight cone}
$\sigma^\vee\subseteq M_{\QQ}$ of $A$ is the cone spanned by all the
lattice vectors $m\in M$ such that $A_m\neq\{0\}$. In the sequel for a
cone $\sigma^\vee\subseteq M_{\QQ}$, we let
$\sigma^\vee_M=\sigma^\vee\cap M$ denote the set of lattice points in
$\sigma^\vee$, so that
$$A=\bigoplus_{m\in \sigma_M^\vee}A_m\chi^{m}\,.$$
Since $A$ is finitely generated, the cone $\sigma^\vee$ is polyhedral
and since the grading is effective, $\sigma^\vee$ is of full dimension
or, equivalently, $\sigma$ is pointed\footnote{A cone in a vector
  space is called pointed if it contains no subspaces of positive
  dimension.}.

An affine $\TT$-variety of complexity 0 is a toric variety. There is a
well known way of describing affine toric varieties in terms of
pointed polyhedral cones in $N_{\QQ}$. To any such cone
$\sigma\subseteq N_{\QQ}$ we associate an affine semigroup algebra
$\KK[\sigma_M^\vee]:=\bigoplus_{m\in\sigma_M^\vee}\KK\chi^m$ and an
affine toric variety $X=\spec\, \KK[\sigma_M^{\vee}]$. Conversely, for
an affine toric variety the corresponding cone $\sigma$ is the dual of
the weight cone $\sigma^\vee$. We note that in this particular case,
$\sigma^\vee\subseteq M_{\QQ}$ is the cone spanned by all lattice
vectors $m\in M$ such that the character $\chi^m:\TT\rightarrow \KK^*$
extends to a regular function on $X$.

In \cite{AlHa}, a combinatorial description of affine $\TT$-varieties
of arbitrary complexity is given. In what follows we recall the main
features of this description specialized to the case of complexity 1
torus actions. In \cite{Ti1} a combinatorial description of complexity
1 actions of reductive groups is given and in \cite{Ti2} it is
specialized for torus actions. For torus actions of complexity 1, the
descriptions in \cite{AlHa} and \cite{Ti1} are equivalent and agree
with the one given earlier (in a slightly more restrictive setting) by
Mumford \cite[Chapters 2 and 4]{KKMS}, cf. \cite{Ti2} and \cite{Vo}.

\begin{definition} \label{stp}
  \begin{enumerate}[(i)]
  \item Let $\sigma$ be a pointed cone in $N_{\QQ}$. We define
    $\pol_{\sigma}(N_{\QQ})$ to be the set of all $\sigma$-tailed
    polyhedra, i.e. polyhedral domains in $N_{\QQ}$ which can be
    decomposed as the Minkowski sum of a compact polyhedron and
    $\sigma$. The set $\pol_{\sigma}(N_{\QQ})$ equipped with the
    Minkowski sum forms an abelian semigroup with neutral element
    $\sigma$.

  \item We let also $\cpl_{\QQ}(\sigma^{\vee})$ denote the set of all
    piecewise linear $\QQ$-valued functions
    $h:\sigma^{\vee}\rightarrow\QQ$ which are upper convex and
    positively homogeneous, i.e.
    $$h(m+m')\geq h(m)+h(m'),\ \mbox{and}\ h(\lambda m)=
    \lambda h(m),\forall m,m'\in \sigma^{\vee},\ \forall\lambda\in
    \QQ_{\geq 0}\,.$$ The set $\cpl_{\QQ}(\sigma^{\vee})$ with the
    usual addition forms an abelian semigroup with neutral element
    $0$.
\end{enumerate}
\end{definition}

For a polyhedron $\Delta\in\pol_{\sigma}(N_{\QQ})$ we define its
support function
\begin{align*} 
  h_{\Delta}:\sigma^{\vee}\rightarrow\QQ,\quad m\mapsto \min\langle
  m,\Delta\rangle\,.
\end{align*}
Clearly, $h_{\Delta}\in \cpl_{\QQ}(\sigma^{\vee})$. The map
$\pol_{\sigma}(N_{\QQ})\rightarrow \cpl_{\QQ}(\sigma^{\vee})$ given by
$\Delta\mapsto h_\Delta$ is an isomorphism of abelian semigroups.

For the following definition we refer to \cite{AlHa}.

\begin{definition} \label{ppd} Let $C$ be a smooth curve. A
  \emph{$\sigma$-polyhedral divisor} on $C$ is a formal sum
  $\DD=\sum_{z\in C}\Delta_z\cdot z$, where
  $\Delta_z\in\pol_{\sigma}(N_{\QQ})$ and $\Delta_z=\sigma$ for all
  but finitely many values of $z$. For $m\in\sigma^{\vee}$ we can
  evaluate $\DD$ in $m$ by letting $\DD(m)$ be the $\QQ$-divisor on
  $C$
  $$\DD(m)=\sum_{z\in C} h_{\Delta_z}(m)\cdot z\,.$$

  A $\sigma$-polyhedral divisor is called \emph{proper} if either $C$
  is affine or $C$ is projective and the following two conditions
  hold.

\begin{enumerate}[(1)]
\item The polyhedron $\deg\DD:=\sum_{z\in C}\Delta_z$ is a proper
  subset of the cone $\sigma$.
\item If $h_{\deg\DD}(m)=0$, then $m$ is contained in the boundary of
  $\sigma^{\vee}$ and a multiple of $\DD(m)$ is principal.
\end{enumerate}
\end{definition}

These two assumptions are counterparts of the conditions that $\DD(m)$
is semiample for all $m\in \sigma^{\vee}_M$ and big for all $m$
contained in the relative interior of $\sigma^{\vee}$,
cf. \cite{AlHa}. They are automatically fulfilled if $C$ is affine.

\begin{definition} \label{nqf} A fan which defines a toric variety
  consists of pointed cones. We need to consider more generally
  objects which we call \emph{quasifans}. These satisfy the usual
  definition of a fan except that their cones are not necessarily
  pointed.

  As usual, for a function $h\in\cpl_{\QQ}(\sigma^{\vee})$ we define
  its \emph{normal quasifan} $\Lambda(h)$ as the coarsest refinement
  of the quasifan of $\sigma^\vee$ such that $h$ is linear in each
  cone $\delta\in\Lambda(h)$. For a $\sigma$-polyhedral divisor $\DD$
  on $C$, we define its normal quasifan $\Lambda(\DD)$ as the coarsest
  common refinement of all $\Lambda(h_{\Delta_z})\,\forall z\in C$. We
  have $\Lambda(\DD)=\Lambda(h_{\deg\DD})$.
\end{definition}

The following theorem gives a combinatorial description of
$\TT$-varieties of complexity 1 analogous to the classical
combinatorial description of toric varieties. This is a specialization
of results in \cite{AlHa} to torus actions of complexity
1. Alternatively, a direct proof is given in \cite{Ti2} for (1) and
(2), while (3) is straightforward from \emph{loc. cit.} See also
Theorem 4.3 in \cite{FlZa1} for the particular case of
$\CC^*$-surfaces.

\begin{theorem} \label{AH-description}
  \begin{enumerate}[(1)]
  \item To any proper $\sigma$-polyhedral divisor $\DD$ on a smooth
    curve $C$ one can associate a normal finitely generated
    effectively $M$-graded domain of dimension $n+1$, where
    $n=\rank(M)$, given by\footnote{For a $\QQ$-divisor $D$, $\lfloor
      D\rfloor$ stands for the integral part and $\{D\}$ for the
      fractional part of $D$.}
    $$A[C,\DD]=\bigoplus_{m\in\sigma^{\vee}_M} A_m\chi^m,\quad
    \mbox{where}\quad A_m=H^0(C,\OO_C(\lfloor \DD(m)\rfloor))\,.$$
  \item Conversely, any normal finitely generated effectively
    $M$-graded domain of dimension $n+1$ is isomorphic to $A[C,\DD]$
    for some smooth curve $C$ and some proper $\sigma$-polyhedral
    divisor $\DD$ on $C$.
  \item Moreover, the $M$-graded domains $A[C,\DD]$ and $A[C',\DD']$
    are isomorphic if and only if $C\simeq C'$, and under this
    identification, $\DD(m)-\DD'(m)$ is linear on $m$ and principal
    for all $m\in\sigma^{\vee}_M$.
  \end{enumerate}
\end{theorem}

In \cite{FiKa} (see also \cite{FlZa1}), all $\CC^*$-surfaces are
divided into three types: elliptic, parabolic and hyperbolic. In the
general case, we will use the following terminology.

An $M$-graded domain $A=A[C,\DD]$ (or, equivalently, a $\TT$-variety
$X$) will be called \emph{elliptic} if $C$ is projective. A
non-elliptic $\TT$-variety will be called \emph{parabolic} if $\sigma$
is of full dimension and \emph{hyperbolic} if $\sigma=\{0\}$. If $\dim
X\geq 3$, this does not exhaust all the possibilities.

\begin{example} \label{ex-hyp} Letting $N=\ZZ^2$ and
  $\sigma=\{(0,0)\}$, in $N_{\QQ}=\QQ^2$ we consider the triangle
  $\Delta_0$ with vertices $(0,0)$,$ (0,1)$ and $(-1/4,-1)$ and the
  segment $\Delta_1=\{0\}\times[0,1]$.
  \begin{figure}[!ht]
    \centering
    \includegraphics[scale=1]{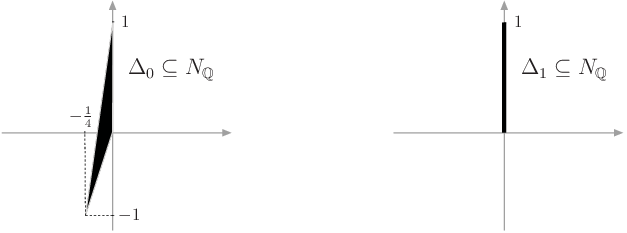}
  \end{figure}

  Let $C=\spec\,\KK[t]$ and
  $\DD=\Delta_0\cdot[0]+\Delta_1\cdot[1]$. In the following picture,
  for the normal quasifans $\Lambda(h_{\Delta_0})$,
  $\Lambda(h_{\Delta_1})$ and $\Lambda(\DD)$ in $M_{\QQ}=\QQ^2$, for
  $i=0,1$ we show the values of $h_i=h_{\Delta_i}$ on each maximal
  cone.
  \begin{figure}[!ht]
    \centering
    \includegraphics[scale=1]{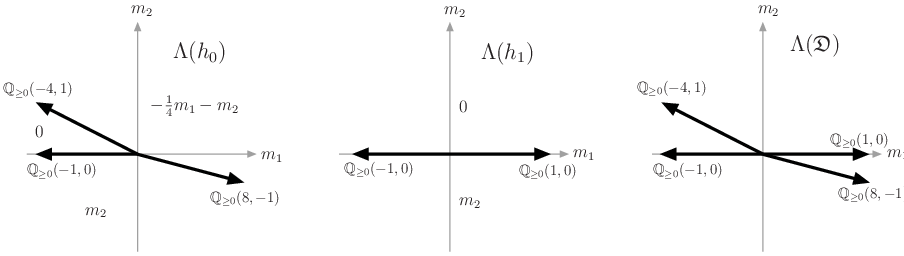}
  \end{figure}

  We let $A=A[C,\DD]$ as in Theorem \ref{AH-description} and
  $X=\spec\,A$. The torus $\TT=(\KK^*)^2$ acts on $X$. Since $C$ is
  affine and $\sigma=\{(0,0)\}$, $X$ is hyperbolic as
  $\TT$-variety. By Theorem \ref{AH-description} we have
  $$A_{(4,0)}=t\KK[t],\quad A_{(-1,0)}=\KK[t],\quad A_{(-4,1)}=\KK[t],
  \quad \mbox{and} \quad A_{(8,-1)}=t(t-1)\KK[t]\,.$$

  An easy calculation shows that the elements
  $$u_1=-t\chi^{(4,0)},\quad u_2=\chi^{(-1,0)},\quad
  u_3=-\chi^{(-4,1)}, \quad \mbox{and} \quad u_4=t(t-1)\chi^{(8,-1)}$$
  generate $A$ as an algebra. Furthermore, they satisfy the
  irreducible relation $u_1+u_1^2u_2^4+u_3u_4=0$, and so
  \begin{align} \label{iso-ex1}
    A\simeq\KK[x_1,x_2,x_3,x_4]/(x_1+x_1^2x_2^4+x_3x_4)\,.
  \end{align}
  The $\ZZ^2$-grading on $A$ is given by $\deg x_1=(4,0)$, $\deg
  x_2=(-1,0)$, $\deg x_3=(-4,1)$, and $\deg x_4=(8,-1)$. The curve $C$
  and the proper polyhedral divisor $\DD$ can be recovered from this
  description of $A$ following the recipe in \cite[Section 11]{AlHa}.

\end{example}

We let $K_0$ denote the function field of $C$. There is a natural
embedding of $M$-graded algebras $A\hookrightarrow K_0[M]$. If $C$ is
affine, then $A_m$ is a locally free $A_0$-module of rank 1 for every
$m\in\sigma^\vee_M$.

Following \cite[Proposition 4.12]{FlZa1}, in the next lemma we show
the way in which our combinatorial description is affected when
passing to a certain cyclic covering.

\begin{lemma} \label{cyclic-divisor} Let $A=A[C,\DD]$, where $C$ is a
  smooth curve with function field $K_0$ and $\DD$ is a proper
  $\sigma$-polyhedral divisor on $C$. Consider the normalization $A'$
  of the cyclic ring extension $A[s\chi^e]$, where $e\in M$, $s^d=f\in
  A_{de}\subseteq K_0$ and $d>0$. Then $A'=A[C',\DD']$, where $C'$ and
  $\DD'$ are defined as follows:
  \begin{enumerate}[(i)]
  \item If $A$ is elliptic, then $A'$ is also elliptic and $C'$ is the
    smooth projective curve with function field $K_0[s]$.

  \item If $A$ is non-elliptic, then $A'$ is also non-elliptic and
    $C=\spec\,A'_0$, where $A'_0$ is the normalization of $A_0$ in
    $K_0[s]$.

  \item In both cases, $\DD'=\sum_{z\in C}\Delta_z\cdot p^*(z),$ where
    $p:C'\rightarrow C$ is the projection.
  \end{enumerate}
\end{lemma}
\begin{proof}
  The normalization $A'$ admits a natural $M$-grading. The latter is
  defined by the $M$-grading on $A$ and by letting $\deg
  s\chi^e=e$. Let $K=\fract\,A$. Since $(s\chi^e)^d-f\chi^{de}=0$,
  $A'$ is the normalization of $A$ in the function field
  $K':=K[s\chi^e]$. But $\chi^{-e}\in K$, so $K'=K[s]$. Moreover
  $K[s]=K_0[s]\otimes \fract\,\KK[M]$, so the function field of $C'$
  is $K_0[s]$, and $A'_0$ is the normalization of $A_0$ in the field
  $K_0[s]$. This proves $(i)$ and $(ii)$.

  For every $m\in M$ we have $\DD'(m)=\sum_{z\in C}
  h_z(m)p^*(z)=p^*(\DD(m))$. Therefore for every $f\in K_0$ there are
  equivalences:
$$\divi_C(f)+\DD(m)\geq0\Leftrightarrow\divi_{C'}(p^*f)+p^*(\DD(m))\geq0\Leftrightarrow \divi_{C'}(f)+\DD'(m)\geq 0\,.$$

Let $m\in\sigma^\vee_M$ and let $r>0$ be such that $\DD(rd\cdot m)$ is
integral. Then
\begin{align*}
  g\in A'_m &\Leftrightarrow g^{rd}\in A_{rdm} \Leftrightarrow \divi_C(g^{rd})+\DD(rd\cdot m)\geq 0 \\
  &\Leftrightarrow \divi_{C'}(g^{rd})+\DD'(rd\cdot m)\geq 0
  \Leftrightarrow \divi_{C'}(g)+\DD'(m)\geq 0\,,
\end{align*}
which proves $(iii)$.
\end{proof}

\subsection{Locally nilpotent derivations and $\KK_+$-actions}

Let $A$ be a commutative ring. A derivation on $A$ is called
\emph{locally nilpotent} (LND for short) if for every $a\in A$ there
exists $n\in\ZZ_{\geq 0}$ such that $\partial^n(a)=0$.

Let $X=\spec\,A$ be an affine variety. Given an LND $\partial$ on $A$,
the map $\phi_\partial:\KK_+\times A\rightarrow A$,
$\phi_\partial(t,f)=e^{t\partial}f$ defines a $\KK_+$-action on $X$,
and any $\KK_+$-action arises in this way. In the following lemma we
collect some well known facts about LNDs over a field of
characteristic 0 not necessarily algebraically closed, needed for
later purposes, see e.g., \cite{Fr2,ML}.
\begin{lemma} \label{LND} Let $A$ be a finitely generated normal
  domain over a field of characteristic 0. If $\partial$ and
  $\partial'$ are two LNDs on $A$, then the following hold:
  \begin{enumerate}[(i)]
  \item $\ker\partial$ is a normal subdomain of codimension 1.
  \item $\ker\partial$ is factorially closed i.e.,
    $ab\in\ker\partial\Rightarrow a,b\in\ker\partial$.
  \item If $a\in A$ is invertible, then $a\in\ker\partial$.
  \item If $\ker\partial=\ker\partial'$, then there exist
    $f,f'\in\ker\partial$ such that $f'\partial=f\partial'$.
  \item For $a\in A$, $\partial a\in(a)\Rightarrow \partial a=0$.
  \item If $a\in\ker\partial$, then $a\partial$ is again a LND.
  \end{enumerate}
\end{lemma}

\begin{definition} \label{LND-equiv} We say that two LNDs $\partial$
  and $\partial'$ on $A$ are \emph{equivalent} if
  $\ker\partial=\ker\partial'$. Geometrically this means that the
  generic orbits of the associated $\KK_+$-actions coincide, cf. also
  Lemma \ref{LND} (iv).
\end{definition}

With dual lattices $M$ and $N$ as in subsection \ref{comb-des}, for a
field extension $\KK\subseteq K_0$ we consider a finitely generated
effectively $M$-graded domain
$$A=\bigoplus_{m\in \sigma_M^\vee}A_m\chi^{m}, \qquad \mbox{where}\qquad A_m\subseteq K_0$$
(we keep our convention from the Introduction regarding $M$-graded
algebras).

A derivation $\partial$ on $A$ is called \emph{homogeneous} if it
sends homogeneous elements into homogeneous elements. Hence $\partial$
sends homogeneous pieces of $A$ into homogeneous pieces.

Let $M_\partial=\{m\in
\sigma^{\vee}_M\mid \partial(A_m\chi^m)\neq0\}$. The action of $\partial$
on homogeneous pieces of $A$ defines a map
$\partial_M:M_{\partial}\rightarrow \sigma^\vee_M$ i.e.,
$\partial(A_m)\subseteq A_{\partial_M(m)}$. By Leibniz rule, for
homogeneous elements $f\in A_m\setminus \ker\partial$ and $g\in
A_{m'}\setminus\ker\partial$ we have
$$\partial(fg)=f\partial(g)+g\partial(f)\in A_{\partial(m+m')},\qquad \partial_M(m+m')=m+\partial_M(m')=m'+\partial_M(m)\,.$$
Thus $\partial_M(m)-m\in M$ is a constant function on
$M_\partial$. This leads to the following definition.

\begin{definition} \label{pn} Let $\partial$ be a nonzero homogeneous
  derivation on $A$. The \emph{degree} of $\partial$ is the lattice
  vector $\deg\partial$ defined by
  $\deg\partial=\deg\partial(f)-\deg(f)$ for any homogeneous element
  $f\notin\ker\partial$.  With this notation the map
  $\partial_M:M_{\partial}\rightarrow \sigma^\vee_M$ is just the
  translation by the vector $\deg\,\partial$.

  We also say that an LND $\partial$ on $A$ is \emph{negative} if
  $\deg\partial\notin \sigma_M^{\vee}$, \emph{non-negative} if
  $\deg\partial\in \sigma_M^{\vee}$, and \emph{positive} if $\partial$
  is non-negative and $\deg\partial\neq 0$.
\end{definition}

It is well known that any LND on $A$ decomposes into a sum of
homogeneous derivations, some of which are locally nilpotent. In lack
of a good reference, in the next lemma we provide a short argument.

\begin{lemma} \label{hLDN} Let $A$ be a finitely generated normal
  $M$-graded domain. For any derivation $\partial$ on $A$ there is a
  decomposition $\partial=\sum_{e\in M}\partial_e$, where $\partial_e$
  is a homogeneous derivation of degree $e$. Moreover, let
  $\Delta(\partial)$ be the convex hull in $M_{\QQ}$ of the set
  $\{e\in M\mid \partial_e\neq 0\}$. Then $\Delta(\partial)$ is a bounded
  polyhedron and for every vertex $e$ of $\Delta(\partial)$,
  $\partial_e$ is locally nilpotent if $\partial$ is.
\end{lemma}
\begin{proof}
  Letting $a_1,\cdots,a_r$ be a set of homogeneous generators of $A$
  we have $A\simeq\KK^{[r]}/I$, where $\KK^{[r]}=\KK[x_1,\cdots,x_r]$
  and $I$ denotes the ideal of relations of $a_1,\cdots,a_r$. The
  $M$-grading and the derivation $\partial$ can be lifted to an
  $M$-grading and a derivation $\partial'$ on $\KK^{[r]}$,
  respectively.

  The proof of Proposition 3.4 in \cite{Fr2} can be applied to an
  $M$-grading, proving that $\partial'=\sum_{e\in M}\partial'_e$,
  where $\partial'_e$ is a homogeneous derivation on
  $\KK^{[r]}$. Furthermore, since $\partial'(I)\subseteq I$ and $I$ is
  homogeneous, we have $\partial'_e(I)\subseteq I$. Hence
  $\partial'_e$ induces a homogeneous derivation $\partial_e$ on $A$
  of degree $e$, proving the first assertion.

  The algebra $A$ being finitely generated, the set $\{e\in
  M\mid\partial_e\neq 0\}$ is finite and so $\Delta(\partial)$ is a
  bounded polyhedron. Let $e$ be a vertex of $\Delta(\partial)$ and
  $n\geq1$. If $ne=\sum_{i=1}^n m_i$ with $m_i\in\Delta(\partial)\cap
  M$, then $m_i=e\ \forall i$. For $a\in A_m\chi^m$ this yields
  $\partial_e^n(a)=\left(\partial^n(a)\right)_{m+ne}$, where
  $\left(\partial^n(a)\right)_{m+ne}$ stands for the summand of degree
  $m+ne$ in the homogeneous decomposition of $\partial^n(a)$. Hence
  $\partial_e$ is locally nilpotent if $\partial$ is so.
\end{proof}

In the following lemma we extend Lemma 1.8 in \cite{FlZa2} to more
general gradings. This lemma shows that any LND $\partial$ on a normal
domain can be extended as an LND to a cyclic ring extension defined by
an element of $\ker\partial$. Actually (i) is contained in
\emph{loc. cit.} while the proof of (ii) is similar and so we omit it.

\begin{lemma} \label{cyclic-LND}
  \begin{enumerate}[(i)]
  \item Let $A$ be a finitely generated normal domain and let
    $\partial$ be an LND on $A$. Given a nonzero element
    $v\in\ker\partial$ and $d>0$, we let $A'$ denote the normalization
    of the cyclic ring extension $A[u]\supseteq A$ in its fraction
    field, where $u^d=v$. Then $\partial$ extends in a unique way to
    an LND $\partial'$ on $A'$.
  \item Moreover, if $A$ is $M$-graded and $\partial$ and $v$ are
    homogeneous, with $\deg v=dm$ for some $m\in M$, then $A'$ is
    $M$-graded as well, and $u$ and $\partial'$ are homogeneous with
    $\deg u=m$ and $\deg\partial'=\deg\partial$.
  \end{enumerate}
\end{lemma}

\begin{sit} \label{18} Recall that $A=\bigoplus_{m\in\sigma^\vee_M}
  A_m\chi^m$, where $A_m\subseteq K_0$, $K_0$ is a field containing
  $\KK$ and $\fract\,A=K_0(M)$ \footnote{For a field $K_0$ and a
    lattice $M$, $K_0(M)$ denotes the function field of
    $K_0[M]$.}. The following lemma provides some useful extension of
  a homogeneous LND $\partial$ on $A$.
\end{sit}

\begin{lemma} \label{LND-tensor} For any homogeneous LND $\partial$ on
  $A$, the following hold:
  \begin{enumerate}[(i)]
  \item The derivation $\partial$ extends in a unique way to a
    homogeneous $\KK$-derivation on $K_0[M]$.
  \item If $\partial(K_0)=0$ then the extension of $\partial$ as in
    (i) restricts to a homogeneous locally nilpotent $K_0$-derivation
    on $K_0[\sigma^\vee_M]$.
  \end{enumerate}
\end{lemma}
\begin{proof}
  The first assertion is evident. Let $\nil(\partial)$ be the
  subalgebra of $K_0[M]$ where $\partial$ acts in a nilpotent way. To
  show (ii), suppose that $\partial(K_0)=0$. Assuming that $f\chi^m\in
  K_0[\sigma^\vee_M]$, we consider $r>0$ such that $A_{rm}\neq
  0$. Letting $g\in A_{rm}$, we have $f^r\chi^{rm}=f'g\chi^{rm}$ for
  some $f'\in K_0$. Thus $f^r\chi^{rm}\in\nil(\partial)$ and so 
  $f\chi^{m}\in\nil(\partial)$.
\end{proof}

In the setting as in the previous lemma, the extension of $\partial$
to $K_0[M]$ will be still denoted by $\partial$. Following
\cite{FlZa2} we use the next definition.

\begin{definition} \label{fh} With $A$ as in \ref{18}, a homogeneous
  LND $\partial$ on A is said to be \emph{of fiber type} if
  $\partial(K_0)=0$ and \emph{of horizontal type} if
  $\partial(K_0)\neq 0$.
\end{definition}

Let $A$ be a finitely generated domain and $X=\spec\,A$. In this
setting, the fact that $\partial$ is homogeneous means that the
corresponding $\KK_+$-action on $X$ is normalized by the $\TT$-action
given by the $M$-grading. Furthermore, $\partial$ is of fiber type if
and only if the general orbits of the corresponding $\KK_+$-action are
contained in the closures of general orbits of the
$\TT$-action. Otherwise, $\partial$ is of horizontal type.

\section{Locally nilpotent derivations on toric varieties}
\label{toric}

In this section we consider more generally toric varieties defined
over a field $\KK$ of characteristic 0, not necessarily algebraically
closed. This will be important in Section \ref{com1} below.

Let $M$ and $N$ be lattices as in Subsection \ref{comb-des}. We also let
$N_{\QQ}=N\otimes\QQ$, $M_{\QQ}=M\otimes\QQ$, and we consider the
natural duality $M_{\QQ}\times N_{\QQ}\rightarrow \QQ$, $(m,p)\mapsto
\langle m,p\rangle$.

\begin{notation} \label{partial-derivatives} Let $\rho\in N$ and $e\in
  M$ be lattice vectors. We define $\partial_{\rho,e}$ as the
  homogeneous derivation of degree $e$ on $\KK[M]$ given by
  $\partial_{\rho,e}(\chi^{m})=\langle m,\rho\rangle\cdot\chi^{m+e}$.
\end{notation}

An easy computation shows that $\partial_{\rho,e}$ is indeed a
derivation. Let $H_\rho$ be the subspace of $M_{\QQ}$ orthogonal to
$\rho$, and $H_\rho^+$ be the halfspace of $M_{\QQ}$ given by
$\langle\cdot,\rho\rangle\geq 0$. The kernel $\ker\partial_{\rho,e}$
is spanned by all characters $\chi^m$ with $m\in M$ orthogonal to
$\rho$, i.e., $\ker\partial_{\rho,e}=\KK[H_\rho\cap M]$.

Let $\nil(\partial_{\rho,e})$ be the subalgebra of $\KK[M]$ where
$\partial_{\rho,e}$ acts in a nilpotent way. Assume that $\langle
e,\rho\rangle=-1$. For every $m\in H_\rho^+\cap M$, the character
$\chi^m\in \nil(\partial_{\rho,e})$ since
$\partial_{\rho,e}^\ell(\chi^m)=0$, where $\ell=\langle
m,\rho\rangle+1$. Thus, the derivation $\partial_{\rho,e}$ restricted
to the subalgebra $\KK[H_\rho^+\cap M]$ is a homogeneous LND. On the
other hand, $\partial_{\rho,e}$ is not locally nilpotent in $\KK[M]$,
in fact for every $m\notin H_\rho^+\cap M$ the character $\chi^m\notin
\nil(\partial_{\rho,e})$ is not nilpotent.

\begin{remark} \label{LND-toric-subalgebra} If $\partial_{\rho,e}$
  stabilizes a subalgebra $A\subseteq\KK[H_\rho^+\cap M]$, then
  $\partial_{\rho,e}|_A$ is also a homogeneous LND on $A$ of degree
  $e$ and $\ker(\partial_{\rho,e}|_A)=A\cap \KK[H_\rho\cap M]$.
\end{remark}

For the rest of this section, we let $\sigma$ be a pointed polyhedral
cone in the vector space $N_{\QQ}$, and
$$A=\KK[\sigma^{\vee}_M]=\bigoplus_{m\in\sigma^{\vee}_M}\KK\chi^m$$
be the affine semigroup algebra of $\sigma$ with the corresponding
affine toric variety $X=\spec\, A$. Since the cone $\sigma$ is
pointed, $\sigma^{\vee}$ is of full dimension and the subalgebra
$A\subseteq \KK[M]$ is effectively graded by $M$.

To every extremal ray $\rho\subseteq\sigma$ we can associate a
codimension 1 face $\tau\subseteq\sigma^\vee$ given by
$\tau=\sigma^{\vee}\cap\rho^\bot$. As usual, we denote an extremal ray
and its primitive vector by the same letter $\rho$. Thus
$\sigma^{\vee}\subseteq H_\rho^+$ and $\tau \subseteq H_\rho$.

\begin{definition}
  Let $\sigma_\rho$ be the cone spanned by all the extremal rays of
  $\sigma$ except $\rho$, so that
  $\sigma^{\vee}=\sigma_\rho^{\vee}\cap H_\rho^+$. We also let
  $$S_{\rho}=\sigma_\rho^\vee\cap\{e\in M\mid\langle e,\rho\rangle=-1\}\,.$$
\end{definition}
This definition can be illustrated on the following picture, where
$\rho\subseteq N_\QQ$ is pointing upwards.
\begin{figure}[!ht]
  \centering
  \includegraphics[scale=1]{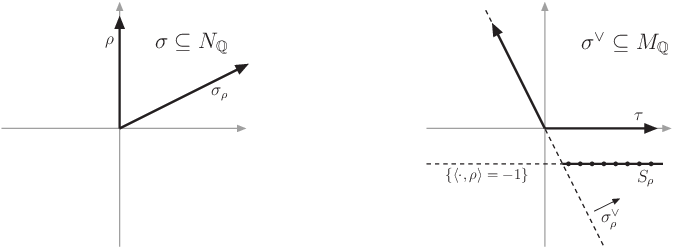}
\end{figure}
\begin{lemma} \label{convex-core} Let $e\in M$. Then $e\in S_\rho$ if
  and only if
  \begin{enumerate}[(i)]
  \item $e\notin\sigma_M^{\vee}$, and
  \item $m+e\in\sigma^\vee_M$, $\forall m\in
    \sigma^{\vee}_M\setminus\tau_M$.
  \end{enumerate}
\end{lemma}

\begin{proof}
  Assume first that $e\in S_\rho$. Then (i) is evident. To show (ii),
  we let $m\in \sigma^{\vee}_M\setminus\tau_M$. Then $m+e\in H_\rho^+$
  because $\langle m+e,\rho\rangle=\langle m,\rho\rangle-1$. Also
  $m\in\sigma^{\vee}\subseteq\sigma_\rho^{\vee}$ yielding
  $m+e\in\sigma_\rho^{\vee}$. Thus
  $m+e\in\sigma^\vee=\sigma_\rho^{\vee}\cap H_\rho^+$.

  To show the converse, we let $e\in M$ be such that (i) and (ii)
  hold. Letting $\rho_i$, $i=1,\cdots,\ell$ be all the extremal rays
  of $\sigma$ except $\rho$, for $m\in\sigma_M^\vee\setminus\tau_M$ we
  have
  $$\langle m+e,\rho_i\rangle=\langle m,\rho_i\rangle+\langle e,\rho_i
  \rangle\geq 0, \ \forall i\in\{1,\cdots,\ell\}\,.$$ If
  $m\in\rho_i^{\bot}\cap\sigma^{\vee}_M$ then $\langle
  m,\rho_i\rangle=0$ and so $\langle e,\rho_i\rangle\geq 0\ \forall
  i$. Thus $e\in\sigma_\rho^{\vee}$. Since
  $e\in\sigma_\rho^{\vee}\setminus\sigma^\vee$, $\langle
  e,\rho\rangle$ is negative. We have $\langle e,\rho\rangle=-1$,
  otherwise $m+e\notin\sigma^\vee$ for any $m\in\sigma^{\vee}_M$ such
  that $\langle m,\rho\rangle=1$. This yields $e\in S_\rho$.
\end{proof}

\begin{remark}\label{toric-non-empty}
  Since $\rho\notin\sigma_\rho$ we have $S_\rho\neq
  \emptyset$. Furthermore, by the previous lemma, $e+m\in S_\rho$
  whenever $e\in S_{\rho}$ and $m\in\tau_M$.
\end{remark}

In the following lemma we provide a translation of Lemma
\ref{convex-core} from the language of convex geometry to that of
affine semigroup algebras.

\begin{lemma} \label{toric-LND} For every pair $(\rho,e)$, where
  $\rho$ is an extremal ray in $\sigma$ and $e$ is a lattice vector in
  $S_\rho$, the homogeneous derivation $\partial_{\rho,e}$ restricts
  to an LND on $A=\KK[\sigma_M^\vee]$ with kernel
  $\ker\partial_{\rho,e}=\KK[\tau_M]$ and $\deg\partial_{\rho,e}=e$.
\end{lemma}

\begin{proof}
  If $\sigma=\{0\}$, then $\sigma$ has no extremal rays, so the
  statement is trivial. We assume in the sequel that $\sigma$ has at
  least one extremal ray $\rho$. By Lemma \ref{convex-core}
  $\partial_{\rho,e}$ stabilizes $A$. Hence by Remark~\ref{LND-toric-subalgebra} (2), $\partial_{\rho,e}$ is a homogeneous
  LND on $A$ with kernel $\KK[\tau_M]$ and of degree $e$.
\end{proof}

The following theorem completes our classification,
cf. \cite[Prop. 11]{De} and \cite[Section 3.4]{Od}.

\begin{theorem} \label{toric-clasification} If $\partial\neq 0$ is a
  homogeneous LND on $A$, then $\partial=\lambda\partial_{\rho,e}$ for
  some extremal ray $\rho$ on $\sigma$, some lattice vector $e\in
  S_\rho$, and some $\lambda\in\KK^*$.
\end{theorem}

\begin{proof}
  The kernel $\ker\partial$ is a subsemigroup subalgebra of $A$ of
  codimension 1. Since $\ker\partial$ is factorially closed (see Lemma
  \ref{LND}), it follows that $\ker\partial=\KK[\sigma^{\vee}_M\cap
  H]$ for a certain codimension 1 subspace $H$ of $M_{\QQ}$.

  If $\sigma^{\vee}\cap H$ is not a codimension 1 face of
  $\sigma^{\vee}$, then $H$ divides the cone $\sigma^\vee$ into two
  pieces. Since the action of $\partial$ on characters is a
  translation by a constant vector $\deg\partial$, only the characters
  from one of these pieces can reach $H$ in a finite number of
  iterations of $\partial$, which contradicts the fact that $\partial$
  is locally nilpotent.

  In the case where $\sigma^{\vee}\cap H=\tau$ is a codimension 1 face
  of $\sigma^{\vee}$, we let $\rho$ be the extremal ray dual to
  $\tau$. Since $\partial$ is an homogeneous LND, the translation by
  $e=\deg\partial$ maps $(\sigma_M^{\vee}\setminus\tau_M)$ into
  $\sigma^\vee_M$. So by Lemma \ref{convex-core}, $e\in S_\rho$ and
  $\partial=\lambda\partial_{\rho,e}$, as required.
\end{proof}

From our classification we obtain the following corollaries.

\begin{corollary} \label{toric-degree} A homogeneous LND $\partial$ on
  a toric variety is uniquely determined, up to a constant factor, by
  its degree.
\end{corollary}
\begin{proof}
  By Theorem \ref{toric-clasification} we have
  $\partial=\lambda\partial_{\rho,e}$ where $e=\deg\partial$. We claim
  that the $\rho$ is uniquely determined by $e$. Indeed, the sets
  $S_{\rho}$ and $S_{\rho'}$ are disjoint for $\rho\neq\rho'$.
  \end{proof}

\begin{corollary} \label{toric-negative} Every homogeneous LND
  $\partial$ on a toric variety $X$ is of fiber type and
  negative\footnote{See Definitions \ref{pn} and \ref{fh}.}.
\end{corollary}
\begin{proof}
  The first claim is evident because $\TT$ acts with an open orbit. By
  Theorem \ref{toric-clasification}, any LND on a toric variety is of
  the form $\lambda\partial_{\rho,e}$. Its degree is
  $\deg\partial_{\rho,e}=e\in S_\rho$ and
  $S_\rho\cap\sigma^\vee=\emptyset$, so $\partial$ is negative.
\end{proof}

\begin{corollary} \label{toric-neq} Two homogeneous LNDs
  $\partial=\lambda\partial_{\rho,e}$ and
  $\partial'=\lambda'\partial_{\rho',e'}$ on $A$ are equivalent if and
  only if $\rho=\rho'$. In particular, there is only a finite number
  of pairwise non-equivalent homogeneous LNDs on $A$.
\end{corollary}
\begin{proof}
  The first assertion follows from the description of
  $\ker\partial_{\rho,e}$ in Lemma \ref{toric-LND} and the second one
  from the fact that $\sigma$, being polyhedral, has only a finite
  number of extremal rays.
\end{proof}

The following corollary shows that the kernel of a homogeneous LND on
a semigroup algebra is finitely generated. Since toric varieties are
rational, this is also a consequence of Theorem 1.2 in \cite{Ku}.

\begin{corollary} \label{fg-toric} Let $X=\spec\ A$ be a toric
  variety. If $\partial:A\rightarrow A$ is a homogeneous LND, then
  $\ker\partial$ is finitely generated as $\KK$-algebra.
\end{corollary}
\begin{proof}
  The corollary follows directly from the description of
  $\ker\partial$ in Lemma \ref{toric-LND}
\end{proof}

\begin{example}
  With $N=\ZZ^3$ we let $\sigma$ be the cone in $N_{\QQ}$ having
  extremal rays $\rho_1=(1,0,0)$, $\rho_2=(0,1,0)$, $\rho_3=(1,0,1)$,
  and $\rho_4=(0,1,1)$. The dual cone $\sigma^\vee\subseteq
  M_{\QQ}=\QQ^3$ is spanned by the lattice vectors $u_1=(1,0,0)$,
  $u_2=(0,1,0)$, $u_3=(0,0,1)$, and $u_4=(1,1,-1)$. Furthermore, these
  elements satisfy the relation $u_1+u_2=u_3+u_4$ and the algebra
  $A=\KK[\sigma_M^\vee]$ is generated by $x_i=\chi^{u_i}$,
  $i=1,\dots,4$. Thus
  \begin{align} \label{iso-toric} A\simeq
    \KK[x_1,x_2,x_3,x_4]/(x_1x_2-x_3x_4)\,.
  \end{align}

  Corollary \ref{toric-neq} shows that there are four non-equivalent
  homogeneous LNDs on $A$ corresponding to the extremal rays
  $\rho_i\subseteq \sigma$. By a routine calculation we obtain
  $$ S_{\rho_1}=\{(-1,b,c)\in M\mid b\geq0, c\geq 1\}, \quad
  S_{\rho_2}=\{(a,-1,c)\in M\mid a\geq0, c\geq 1\},$$
  $$S_{\rho_3}=\{(a,b,c)\in M\mid a\geq0, b+c\geq 0, a+c=-1\},\mbox{ and}$$
  $$S_{\rho_4}=\{(a,b,c)\in M\mid b\geq0, a+c\geq 0, b+c=-1\}\,.$$

  Letting $e_1=(-1,0,1)$, $e_2=(0,-1,1)$, $e_2=(0,1,-1)$,
  $e_4=(1,0,-1)$, $\partial_i=\partial_{\rho_i,e_i}$, and
  $m=(m_1,m_2,m_3)$, we have
  $$\partial_1(\chi^m)=m_1\cdot\chi^{m+e_1},\quad 
  \partial_2(\chi^m)=m_2\cdot\chi^{m+e_2},$$
  $$\partial_3(\chi^m)=(m_1+m_3)\cdot\chi^{m+e_3},\quad
  \mbox{and} \quad
  \partial_4(\chi^m)=(m_2+m_3)\cdot\chi^{m+e_4}\,,
  $$
  Finally, under the isomorphism of \eqref{iso-toric} the four
  homogeneous LNDs on $A$ are given by
  $$\partial_1=x_3\frac{\partial}{\partial x_1}+x_2\frac{\partial}{\partial x_4}, \quad
  \partial_2=x_3\frac{\partial}{\partial
    x_2}+x_1\frac{\partial}{\partial x_4},
  $$
  $$\partial_3=x_4\frac{\partial}{\partial
    x_1}+x_2\frac{\partial}{\partial x_3}, \quad \mbox{and}\quad
  \partial_4=x_4\frac{\partial}{\partial
    x_2}+x_1\frac{\partial}{\partial x_3}\,.
  $$
\end{example}

\section{Locally nilpotent derivations on $\TT$-varieties of
  complexity 1}
\label{com1}

In this section we give a complete classification of homogeneous LNDs
on $\TT$-varieties of complexity 1 over an algebraically closed field
$k$ of characteristic 0. In the first part we treat the case of a
homogeneous LNDs of fiber type, while in the second one we deal with
the more delicate case of homogeneous LNDs of horizontal type.

We fix the $n$-dimensional torus $\TT$, a smooth curve $C$ and a
proper $\sigma$-polyhedral divisor $\DD=\sum_{z\in C} \Delta_z\cdot z$
on $C$. Letting $K_0$ be the function field of $C$, we consider the
affine variety $X=\spec\,A$, where
$$A=A[C,\DD]=\bigoplus_{m\in\sigma^{\vee}_M}A_m\chi^m,\quad\mbox{with}\quad A_m=H^0\left(C,\OO(\lfloor\DD(m)\rfloor)\right)\subseteq K_0\,.$$
We denote by $h_z=h_{\Delta_z}$ the support function of $\Delta_z$ so
that $\DD(m)=\sum_{z\in C} h_z(m)\cdot z$. We also fix a homogeneous
LND $\partial$ on $A$.

In this context, we can interpret Definitions \ref{pn} and \ref{fh} as
follows.

\begin{lemma} \label{fiber-negative} With the notation as above, let
  $\partial$ be a homogeneous LND on $A$. Then the following hold.
  \begin{enumerate}[(i)]
  \item If $\partial$ is of fiber type, then $\partial$ is negative
    and $\ker\partial=\bigoplus_{m\in \tau_M}A_m\chi^m$, where $\tau$
    is a codimension 1 face of $\sigma^{\vee}$.
  \item Assuming further that $A$ is non-elliptic, $\partial$ is of
    fiber type if and only if $\partial$ is negative.
  \end{enumerate}
\end{lemma}

\begin{proof}
  To prove (i) we let $\partial$ be a homogeneous LND of fiber type on
  $A$. By Lemma \ref{LND-tensor} we can extend $\partial$ to a
  homogeneous LND $\bar{\partial}$ on $\bar{A}=K_0[\sigma^\vee_M]$
  which is an affine semigroup algebra over $K_0$. Since
  $\partial(K_0)=0$, $\bar{\partial}$ is a locally nilpotent
  $K_0$-derivation. It follows from Corollary \ref{toric-negative}
  that $\deg\partial=\deg\bar{\partial}\notin\sigma_M^\vee$, so
  $\partial$ is negative.

  Furthermore, Lemma \ref{toric-LND} and Theorem
  \ref{toric-clasification} show that
  $\ker\bar{\partial}=K_0[\tau_M]$, where $\tau$ is a codimension 1
  face of $\sigma^{\vee}$. Thus
$$\ker\partial=A\cap\ker\bar{\partial}=\bigoplus_{m\in\tau_M}(A_m\cap K_0)\chi^m=\bigoplus_{m\in\tau_M}A_m\chi^m\,,$$
which proves (i).

To prove (ii) we assume further that $A$ is non-elliptic. Let
$\partial$ be a negative homogeneous LND on $A$. Let $\bar{\partial}$
be the extension of $\partial$ to $K_0[M]$ as in Lemma
\ref{LND-tensor}. Since $\partial$ is negative,
$\partial(A_0)\subseteq A_{\deg\partial}=0$. Since $A$ is non-elliptic
we have $K_0=\fract\, A_0$, so $\bar\partial(K_0)=0$ and $\partial$ is
of fiber type.
\end{proof}

\begin{remark}
  In the elliptic case, the second assertion in Lemma
  \ref{fiber-negative} does not hold, in general. Consider for
  instance the elliptic $\KK$-domain $A=\KK[x,y]$ graded via $\deg
  x=\deg y=1$. Then the partial derivative $\partial_x$ is a negative
  homogeneous LND of horizontal type on $A$.
\end{remark}

\subsection{Homogeneous LNDs of fiber type}

In this subsection we consider a homogeneous LND $\partial$ on $A$ of
fiber type. Let as before $\bar{A}=K_0[\sigma^\vee_M]$ be the affine
semigroup $K_0$-algebra with cone $\sigma\in N_{\QQ}$ over the field
$K_0$ of rational functions of $C$. By Lemma \ref{LND-tensor},
$\partial$ can be extended to a homogeneous locally nilpotent
$K_0$-derivation on $\bar{A}$. To classify homogeneous LNDs of fiber
type, we will rely on the classification of homogeneous LNDs on affine
semigroup algebras from the previous section.

If $\sigma$ has no extremal ray then $\sigma=0$ and
$\sigma^\vee=M_{\QQ}$. By Lemma \ref{fiber-negative} in this case
there are no homogeneous LND of fiber type. So we may assume in the
sequel that $\sigma$ has at least one extremal ray, say $\rho$. Let
$\tau$ be its dual codimension 1 face, and let $S_\rho$ be as defined
in Lemma \ref{convex-core}.
\begin{lemma} \label{3.4} For any $e\in S_\rho$,
$$D_e:=\sum_{z\in C}\max_{m\in\sigma^\vee_M\setminus\tau_M}(h_z(m)-h_z(m+e))\cdot z$$
is a well defined $\QQ$-divisor on $C$.
\end{lemma}

\begin{proof} 
  By Lemma \ref{convex-core}, for all $m\in
  \sigma^\vee_M\setminus\tau_M$, $m+e$ is contained in $\sigma^\vee_M$
  and thus $h_z(m)$ and $h_z(m+e)$ are well defined. Recall that for
  any $z\in C$, the function $h_z$ is upper convex and piecewise
  linear on $\sigma^\vee$. Thus the above maximum is achieved by one
  of the linear pieces of $h_z$ i.e., by one of the maximal cones in
  the normal quasifan $\Lambda(h_z)$ (see Definition \ref{nqf}).

  For every $z\in C$, we let
  $\{\delta_{1,z},\cdots,\delta_{\ell_z,z}\}$ be the set of all
  maximal cones in $\Lambda(h_z)$ and $g_{r,z},\
  r\in\{1,\cdots,\ell_z\}$ be the linear extension of
  $h_z|_{\delta_{r,z}}$ to $M_{\QQ}$. Since the maximum is achieved by
  one of the linear pieces we have
$$\max_{m\in\sigma^\vee_M\setminus\tau_M}(h_z(m)-h_z(m+e))=\max_{r\in\{1,\cdots,\ell_z\}}(-g_{r,z}(e))\,.$$
Since $g_{r,z}(e)\in\QQ\ \forall (r,z)$, $D_e$ is indeed a
$\QQ$-divisor.
\end{proof}

\begin{remark} \label{De-simple} With the notation as in the preceding
  proof we can provide a better description of $D_e$. Since $\tau$ is
  a codimension 1 face of $\sigma^\vee$, it is contained as a face in
  one and only one maximal cone $\delta_{r,z}$. We may assume that
  $\tau\subseteq\delta_{1,z}$.  By the upper convexity of $h_z$ we
  have $\ g_{1,z}(e)\leq g_{r,z}(e)$ $\forall r$ and so
  $D_e=-\sum_{z\in C} g_{1,z}(e)\cdot z$.
\end{remark}

\begin{notation} \label{3.6} We let
  $\Phi_e=H^0(C,\OO_C(\lfloor-D_e\rfloor))$. Thus for any
  $\varphi\in\Phi_e$ and any $m\in\sigma^\vee_M\setminus\tau_M$ we
  have $$\divi(\varphi)\geq \left\lceil D_e\right\rceil\geq D_e\geq
  \sum_{z\in C}(h_z(m)-h_z(m+e))\cdot z=\DD(m)-\DD(m+e)\,.$$
\end{notation}

There is a natural way to associate to a nonzero function
$\varphi\in\Phi_e$ a homogeneous LND of fiber type on $A$. More
precisely we have the following lemma.

\begin{lemma} \label{fib-core} To any triple $(\rho,e,\varphi)$, where
  $\rho$ is an extremal ray of $\sigma$, $e\in S_\rho$ is a lattice
  vector, and $\varphi\in\Phi_e$ is a nonzero function, we can
  associate a homogeneous LND $\partial_{\rho,e,\varphi}$ on
  $A=A[C,\DD]$ with kernel
$$\ker\partial_{\rho,e,\varphi}=\bigoplus_{m\in\tau_M}A_m\chi^m,\qquad\mbox{and}\qquad \deg\partial_{\rho,e,\varphi}=e\,.$$
\end{lemma}

\begin{proof}
  Letting $\bar{A}=K_0[\sigma^\vee_M]$, we consider the $K_0$-LND
  $\partial_{\rho,e}$ on $\bar{A}$ as in Lemma \ref{toric-LND}. Since
  $\varphi\in K_0$, $\varphi\partial_{\rho,e}$ is again an $K_0$-LND
  on $\bar{A}$.

  We claim that $\varphi\partial_{\rho,e}$ stabilizes $A\subseteq
  \bar{A}$. Indeed, let $f\in A_m\subseteq K_0$ be a homogeneous
  element so that $\divi\,f+\lfloor\DD(m)\rfloor\geq 0$. If
  $m\in\tau_M$, then $\varphi\partial_{\rho,e}(f\chi^m)=0$. If
  $m\in\sigma^\vee_M\setminus\tau_M$, then
  $$\varphi\partial_{\rho,e}(f\chi^m)=\varphi
  f\partial_{\rho,e}(\chi^m)= m_0\varphi f\chi^{m+e}\,,$$ where
  $m_0:=\langle m,\rho\rangle\in \ZZ_{>0}$. Moreover by virtue of
  Notation \ref{3.6},
\begin{align*}
  \divi(m_0\varphi f)+\lfloor\DD(m+e)\rfloor&=\divi\,\varphi+\divi\,f+\lfloor\DD(m+e)\rfloor \\
  &\geq\DD(m)-\DD(m+e)-\lfloor\DD(m)\rfloor+\lfloor\DD(m+e)\rfloor \\
  &=\{\DD(m)\}-\{\DD(m+e)\}\,.
\end{align*}
Since the divisor $\divi(m_0\varphi f)+\lfloor\DD(m+e)\rfloor$ is
integral and all the values of the divisor $\{\DD(m)\}-\{\DD(m+e)\}$
are in the interval $]-1,1[$ we have
$$\divi(m_0\varphi f)+\lfloor\DD(m+e)\rfloor\geq 0\quad \mbox{ and so } \quad m_0\varphi f\in A_{m+e}\,,$$
yielding the claim. Finally
$\partial_{\rho,e,\varphi}:=\varphi\partial_{\rho,e}|_A$ is an
homogeneous LND on $A$ with kernel
$\ker\partial_{\rho,e,\varphi}=\bigoplus_{m\in\tau_M}A_m\chi^m$, as
desired.
\end{proof}

\begin{remark} \label{fiber-phi} We have shown actually that for every
  $\varphi\in \Phi_e$, $\varphi A_m\subseteq A_{m+e}$ for all
  $m\in\sigma^\vee_M\setminus\tau_M$. It is easily seen from the
  construction of the divisor $D_e$ that all the functions $\varphi\in
  K_0$ satisfying this property are contained in $\Phi_e$.
\end{remark}

The following theorem gives the converse of Lemma \ref{fib-core} and
so completes our classification of homogeneous LNDs of fiber type on
$\TT$-varieties.

\begin{theorem} \label{fiber-clas} Every nonzero homogeneous LND
  $\partial$ of fiber type on $A=A[C,\DD]$ is of the form
  $\partial=\partial_{\rho,e,\varphi}$ for some extremal ray
  $\rho\subseteq\sigma$, some lattice vector $e\in S_\rho$, and some
  function $\varphi\in\Phi_e$.
\end{theorem}

\begin{proof}
  Since $\partial$ is of fiber type, $\partial|_{K_0}=0$ and so
  $\partial$ can be extended to a $K_0$-LND $\bar{\partial}$ on the
  affine semigroup algebra $\bar{A}=K_0[\sigma^\vee_M]$. By Theorem
  \ref{toric-clasification} we have
  $\bar{\partial}=\varphi\partial_{\rho,e}$ for some extremal ray
  $\rho$ of $\sigma$, some $e\in S_\rho$ and some $\varphi\in
  K_0$. Since $A$ is stable under $\varphi\partial_{\rho,e}$, by
  Remark \ref{fiber-phi}, $\varphi\in\Phi_e$ and so
  $\partial=\varphi\partial_{\rho,e}|_A=\partial_{\rho,e,\varphi}$.
\end{proof}

\begin{corollary}
  Let as before $X=\spec\,A$ be a $\TT$-variety of complexity 1,
  $\partial$ be a homogeneous LND of fiber type on $A$, and let
  $f\chi^m\in A\setminus\ker\partial$ be a homogeneous element. Then
  $\partial$ is completely determined by the image
  $g\chi^{m+e}:=\partial(f\chi^m)\in A_{m+e}\chi^{m+e}$.
\end{corollary}
\begin{proof}
  By the previous theorem $\partial=\partial_{\rho,e,\varphi}$ for
  some extremal ray $\rho$, some $e\in S_\rho$, and some
  $\varphi\in\Phi_e$, where $e=\deg\partial$ and $\rho$ is uniquely
  determined by $e$, see Corollary \ref{toric-degree}.

  In the proof of Lemma \ref{fib-core} it was shown that
  $\partial_{\rho,e,\varphi}(f\chi^m)=m_0\varphi f\chi^{m+e}$. Thus
  $\varphi=\frac{g}{m_0f}\in K_0$ is also uniquely determined by our
  data.
\end{proof}

\begin{corollary} \label{feq} Two homogeneous LND
  $\partial=\partial_{\rho,e,\varphi}$ and
  $\partial'=\partial_{\rho',e',\varphi'}$ of fiber type on $A$ are
  equivalent if and only if $\rho=\rho'$. In particular, there is a
  finite number of pairwise non-equivalent LNDs of fiber type on $A$.
\end{corollary}
\begin{proof}
  The first assertion follows from the description of
  $\ker\partial_{\rho,e,\varphi}$ in Lemma \ref{fib-core}. The second
  one follows from the fact that $\sigma$ has a finite number of
  extremal rays.
\end{proof}

In the following proposition we show that the kernel of a homogeneous
LND of fiber type is finitely generated.

\begin{proposition} \label{fg-fib} Let $\partial$ be a homogeneous LND
  on $A=A[C,\DD]$, where $\DD$ is a proper polyhedral $\sigma$-divisor
  on a smooth curve $C$. If $\partial$ is of fiber type, then
  $\ker\partial$ is finitely generated.
\end{proposition}

\begin{proof}
  In the notation of Theorem \ref{fiber-clas} we have
  $\partial=\partial_{\rho,e,\varphi}$, where $\rho\subseteq\sigma$ is
  an extremal ray. Letting $\tau\subseteq\sigma^\vee$ be the
  codimension 1 face dual to $\rho$, Lemma \ref{fib-core} shows that
  $\ker\partial=\bigoplus_{m\in\tau_M} A_m\chi^m$.

  Let $a_1,\ldots,a_r$ be a set of homogeneous generators of
  $A$. Without loss of generality, we assume further that $\deg
  a_i\in\tau_M$ if and only if $1\leq i\leq s<r$. We claim that
  $a_1,\ldots,a_s$ generate $\ker\partial$. Indeed, let $P$ be any
  polynomial such that $P(a_1,\ldots,a_r)\in\ker\partial$. Since
  $\tau\subseteq\sigma^\vee$ is a face, $\sum m_i\in \tau_M$ for
  $m_i\in\sigma_M^\vee$ implies that $m_i\in\tau$ $\forall i$. Hence
  all the monomials composing $P(a_1,\ldots,a_r)$ are monomials in
  $a_1,\ldots,a_s$, proving the claim.
\end{proof}

Given an extremal ray $\rho\subseteq \sigma$ and $e\in S_\rho$, it
might happen that $\dim \Phi_e=h^0(C,\OO_C(\lfloor -D_e\rfloor))=0$,
so that there exist no homogeneous LND $\partial$ of fiber type on $A$
with $\deg\partial=e$ and
$\ker\partial=\bigoplus_{m\in\tau_M}A_m\chi^m$. In the following lemma
we give a criterion for the existence of $e\in S_\rho$ such that
$\dim\Phi_e$ is nonzero.

\begin{lemma} \label{conv-fib} Let $A=A[C,\DD]$, and let
  $\rho\subseteq\sigma$ be an extremal ray dual to a codimension one
  face $\tau\subseteq \sigma^\vee$. There exists $e\in S_\rho$ such
  that $\dim\Phi_e$ is positive if and only if the curve $C$ is affine
  or $C$ is projective and $h_{\deg\DD}|_{\tau}\not\equiv 0$.
\end{lemma}
\begin{proof}
  If $C$ is affine, then for any $\ZZ$-divisor $D$ the sheaf
  $\OO_C(D)$ is generated by the global sections. It follows in this
  case that $\dim\Phi_e>0$.

  Let further $C$ be a projective curve of genus $g$. If $\deg
  \lfloor-D_e\rfloor<0$ then $\dim\Phi_e=0$. On the other hand, by the
  Riemann-Roch theorem $\dim\Phi_e>0$ if $\deg \lfloor-D_e\rfloor\geq
  g$ (see Lemma 1.2 in \cite[Chapter IV]{Ha}).

  Letting $h=h_{\deg\DD}=\sum_{z\in C}h_z$, with the notation of
  Remark \ref{De-simple} we have $h|_{\tau}=\sum_{z\in C} g_{1,z}$ and
  $\deg(-D_e)=\sum_{z\in C}g_{1,z}(e)$. By the definition of proper
  $\sigma$-polyhedral divisor, $h(m)>0$ for any $m$ in the relative
  interior of $\sigma^\vee$.

  If $h|_{\tau}\equiv0$ then by the linearity of $g_{1,z}$ we obtain
  that $\deg(-D_e)<0$, so $\deg \lfloor-D_e\rfloor < 0$ and
  $\dim\Phi_e=0$.

  If $h|_{\tau}\not\equiv0$ then by the upper convexity of $h$,
  $h(m)>0$ for all $m$ in the relative interior of $\tau$. By Remark
  \ref{De-simple}, $\deg(-D_e)$ is linear on $e$ and so, according to
  Remark \ref{toric-non-empty}, we can choose a suitable $e\in S_\rho$
  so that $\deg \lfloor-D_e\rfloor \geq g$. Hence $\dim\Phi_e>0$.
\end{proof}

We can now deduce the following corollary.

\begin{corollary} \label{fiex} Let $A=A[C,\DD]$, and let
  $\rho\subseteq\sigma$ be an extremal ray dual to a codimension one
  face $\tau\subseteq \sigma^\vee$. There exists a homogeneous LND of
  fiber type $\partial$ on $A$ such that
  $\ker\partial=\bigoplus_{m\in\tau_M}A_m\chi^m$ if and only if $C$ is
  affine or $C$ is projective and $\rho\cap\deg\DD=\emptyset$.
\end{corollary}
\begin{proof}
  Since $\rho\cap\deg\DD=\emptyset$ is equivalent to
  $h_{\deg\DD}|_{\tau}\not\equiv 0$, the corollary follows from
  Theorem \ref{fiber-clas} and Lemma \ref{conv-fib}.
\end{proof}

\begin{remark} \label{feq2} By Corollaries \ref{feq} and \ref{fiex},
  the equivalence classes of LNDs of fiber type on $A=A[C,\DD]$ are in
  one to one correspondence with the extremal rays $\rho\subseteq
  \sigma$ if $C$ is affine and with extremal rays $\rho\subseteq
  \sigma$ such that $\rho\cap\deg\DD=\emptyset$ if $C$ is projective.
\end{remark}

\begin{remark}
  In the recent preprint \cite{Li1}, we generalize the methods of this
  section to give a classification of LNDs of fiber type in arbitrary
  complexity.
\end{remark}

\subsection{Homogeneous LNDs of horizontal type}

Let $A=A[C,\DD]$, where $\DD$ is a proper $\sigma$-polyhedral divisor
on a smooth curve $C$. We consider a homogeneous LND $\partial$ of
horizontal type on $A$. We also denote by $\partial$ its extension to
a homogeneous $\KK$-derivation on $K_0[M]$, where $K_0$ is the field
of rational functions of $C$ (see Lemma \ref{LND-tensor} (i)).

The existence of a homogeneous LND of horizontal type imposes strong
restrictions on $C$, as we show in the next lemma.

\begin{lemma} \label{hor-freeness} If there exists a homogeneous LND
  $\partial$ of horizontal type on $A=A[C,\DD]$, then $C\simeq\PP^1$
  in the case where $A$ is elliptic and $C\simeq\AF^1$ in the case
  where $A$ is non-elliptic. In the latter case $A_m$ is a free
  $A_0$-module of rank 1 for every $m\in\sigma^\vee_M$ and so
$$A_m=\varphi_mA_0\quad \mbox{for some} \quad \varphi_m\in A_m \quad \mbox{such that} \quad \divi(\varphi_m)+\lfloor\DD(m)\rfloor=0\,.$$
\end{lemma}

\begin{proof}
  Let $\pi:X=\spec\,A\dashrightarrow C$ be the rational quotient for
  the $\TT$-action given by the inclusion $\pi^*:K_0\hookrightarrow
  K=\fract\,A$. Since $X$ is normal, the indeterminacy locus $X_0$ of
  $\pi$ has codimension greater than 1, and so the general orbits of
  the $\KK_+$-action corresponding to $\partial$ are contained in
  $X\setminus X_0$.

  Since $\partial|_{K_0}\neq 0$, the general orbits of the
  $\KK_+$-action on $X$ are not contained in the fibers of $\pi$, so
  map dominantly onto $C$. Hence $C$ being dominated by $\AF^1$ we
  have $C\simeq\PP^1$ in the elliptic case and $C\simeq\AF^1$ in the
  non-elliptic case.

  Thus, if $C$ is affine then $A_0=\KK[t]$ and so $A_m$ is a locally
  free $A_0$-module of rank 1 for any $m\in\sigma^\vee_M$. By the
  primary decomposition, any locally free module over a principal ring
  is free and so $A_m\simeq A_0$ as a module (see also Ch. VII \S4
  Corollary 2 in \cite{Bu}). Now the last assertion easily follows.
\end{proof}

\begin{sit}
  For the rest of this section we let $K_0=\KK(t)$, $C=\PP^1$ in the
  elliptic case, and $C=\AF^1$ otherwise. We also let $S_\partial$ be
  the set of all lattice vectors $m\in M$ such that $\ker\partial\cap
  A_m\chi^m\neq\{0\}$, $L(\partial)\subseteq M$ be the sublattice
  spanned by $S_\partial$, and $\omega^{\vee}(\partial)$ be the cone
  spanned by $S_\partial$ in $M_{\QQ}$. We write $L$ and $\omega^\vee$
  instead of $L(\partial)$ and $\omega^\vee(\partial)$ whenever
  $\partial$ is clear from the context.
\end{sit}

\begin{lemma} \label{kernel-linear} Let $A=A[C,\DD]$, where $\DD$ is a
  proper $\sigma$-polyhedral divisor on $C$, and let $\partial$ be a
  homogeneous LND of horizontal type on $A$. With the notation as
  above, the following hold.
  \begin{enumerate}[(1)]
  \item The kernel $\ker\partial$ is a semigroup algebra given by
    $\ker\partial=\bigoplus_{m\in\omega^{\vee}_L} \KK\varphi_m\chi^m$,
    where $\varphi_m\in A_m$.
  \item For all $m\in\omega^{\vee}_L$, in the non-elliptic case
    $\divi(\varphi_m)+\DD(m)=0$, while in the elliptic one
    $\divi(\varphi_m)+\DD(m)=\lambda\cdot[z_\infty]$ for some
    $z_\infty\in\PP^1$ and some positive $\lambda\in\QQ$.
  \item The cone $\omega^{\vee}$ is a maximal cone of the quasifan
    $\Lambda(\DD)$ in the non-elliptic case, and of the quasifan
    $\Lambda(\DD|_{\PP^1\setminus\{z_\infty\}})$ in the elliptic
    one. In particular, $\rank(L)=n$.
  \item $M$ is spanned by $\deg\partial$ and $L$. More precisely, any
    $m\in M$ can be uniquely written as $m=l+r\deg\partial$ for some
    $l\in L$ and some $r\in\ZZ$ with $0\leq r<d$, where $d>0$ is the
    smallest integer such that $d\deg\partial\in L$.
  \end{enumerate}
\end{lemma}

\begin{proof}
  Since $\KK\subseteq\ker\partial$ we have $0\in S_\partial$. If
  $m,m'\in S_\partial$ then $m+m'\in S_\partial$ and so $S_\partial$
  is a subsemigroup of $\sigma^{\vee}_M$.
 
  For any $f\in K_0=\KK(t)$ we have $\partial(f)=f'(t)\partial(t)$,
  where $\partial(t)\neq 0$ since $\partial$ is of horizontal
  type. Thus $\partial(f)=0$ if and only if $f$ is constant. Let us
  fix $m\in S_{\partial}$. If $\varphi_m,\varphi_m'\in
  \ker\partial\cap A_m\chi^m$ are nonzero, then
  $\varphi_m/\varphi_m'\in \ker\partial\cap K_0=\KK$ and so
  $\varphi_m'=\lambda \varphi_m$ for some $\lambda\in\KK^*$.

  Hence $\ker\partial=\bigoplus_{m\in S_{\partial}}\KK\varphi_m\chi^m$
  and $\ker\partial$ is a semigroup algebra. Since $\ker\partial$ is
  normal, $S_\partial$ is saturated, and so
  $S_\partial=\omega^{\vee}_L$, which proves (1).

  To prove (2), we assume first that $C$ is affine. Given
  $m\in\omega_L^\vee$, we let $\varphi_m$ be as in Lemma
  \ref{hor-freeness}. Since $\ker\partial$ is factorially closed, if
  $f\varphi_m\chi^m\in\ker\partial\cap A_m\chi^m$ for some $f\in A_0$,
  then $f\in\ker\partial\cap A_0=\KK$ and
  $\varphi_m\chi^m\in\ker\partial\cap A_m\chi^m$. The latter implies
  that $\varphi_m^r\chi^{rm}\in\ker\partial\cap A_{rm}\chi^{rm}$
  $\forall r\geq 1$, and so $r\lfloor \DD(m)\rfloor=\lfloor
  r\DD(m)\rfloor$ $\forall r\geq 1$. Hence $\DD(m)$ is an integral
  divisor, which yields (2) in the non-elliptic case.

  In the case where $C=\PP^1$, we may suppose that that
  $z_\infty=\infty$. Given $m\in \omega_L^\vee$, let us assume that
  $\divi(\varphi_m)+\lfloor\DD(m)\rfloor\geq [0]+[\infty]$ so that
  $t\varphi_m\in A_m$ and $t^{-1}\varphi_m\in A_m$. We have
  $(t\varphi_m\chi^m)(t^{-1}\varphi_m\chi^m)=(\varphi_m\chi^m)^2\in\ker\partial$. Thus
  $t\varphi_m\chi^m\in\ker\partial$, which contradicts (1). Henceforth
  $\divi(\varphi_m)+\lfloor\DD(m)\rfloor=\lambda\cdot[z_\infty]$,
  $\lambda\in\ZZ_{\geq0}$. An argument similar to that employed in the
  non-elliptic case, yields
  $\divi(\varphi_m)+\DD(m)=\lambda\cdot[z_\infty]$ for some positive
  $\lambda\in\QQ$, proving (2).

  We have $\dim\ker\partial=\dim\omega^{\vee}$. Since $\partial$ is an
  LND, $\ker\partial$ has codimension 1 in $A$. Hence $\omega^{\vee}$
  is of full dimension in $M_{\QQ}$. Furthermore, in the non-elliptic
  case (2) shows that $h_z|_{\omega^\vee}$ is linear $\forall
  z\in\AF^1$, so that $\omega^\vee$ is contained in a maximal cone
  $\delta$ in $\Lambda(\DD)$.

  Assume that $\omega^\vee\subsetneq\delta$. Let
  $m\in\delta\setminus\omega^\vee$ and $\varphi_m\in\KK(t)$ be such
  that $\DD(m)$ is integral and $\divi(\varphi_m)+\DD(m)=0$. Letting
  $m'\in\omega^\vee_L$ be such that $m+m'\in\omega^\vee_L$, the
  linearity of $\DD$ implies
  $\varphi_m\chi^m\varphi_{m'}\chi^{m'}=\varphi_{m+m'}\chi^{m+m'}\in\ker\partial$. Hence
  $\varphi_m\chi^m\in\ker\partial$ which is a contradiction, proving
  (3) in the non-elliptic case. In the elliptic case a similar
  argument (with $z\in\PP^1\setminus\{z_\infty\}$) provides the
  result.

  Finally, since $\sigma_M^\vee$ spans $M$ as a lattice and $\partial$
  is a homogeneous LND, for any $m\in M$ we have $m+r\deg\partial\in
  L$ for some $r\in\ZZ$. Thus for $0\geq r>-d$ the decomposition as in
  (4) is unique because of the minimality of $d$.
\end{proof}

The following corollary shows that the kernel of a homogeneous LND on
$A$ is a semigroup algebra and so the kernel is finitely
generated. Since, by Lemma \ref{hor-freeness}, $\spec\,A$ is rational,
this is also a consequence of Theorem 1.2 in \cite{Ku}.

\begin{corollary} \label{fgraro} In the notation of Lemma
  \ref{kernel-linear}, by (3) $\omega\subseteq N_\QQ$ is a pointed
  polyhedral cone and by (1)
$$\ker\partial=\bigoplus_{m\in\omega^{\vee}_L} \KK\varphi_m\chi^m\simeq \KK[\omega_L^\vee]$$
is an affine semigroup algebra, in particular $\ker\partial$ is
finitely generated.
\end{corollary}

Let us consider two basic examples, one with a non-elliptic
$\TT$-action and the other one with an elliptic $\TT$-action. They are
universal in the sense of Lemma \ref{hor-ome} below. We use both
examples in our final classification, cf. Lemma \ref{pres} and Theorem
\ref{resultado}.

Starting with an affine toric variety $X$ and a homogeneous LND
$\partial$ of fiber type (see Corollary \ref{toric-negative}), we can
restrict the big torus action to an appropriate codimension 1 subtorus
$\TT$ so that $\partial$ becomes of horizontal type for the
$\TT$-action of complexity 1 on $X$. This is actually the case in our
examples.

\begin{example} \label{hor-fund-exa} Letting $A=A[C,\DD]$, where
  $C=\AF^1$, $p\in N_{\QQ}$, and $\DD=(p+\sigma)\cdot [0]$ we have
  that $h_0:\sigma^\vee\rightarrow\QQ,\, m\mapsto\langle m,p\rangle$
  is linear and $h_z=0$ $\forall z\in\KK^*$. Denoting by
  $h:M_{\QQ}\rightarrow \QQ$ the linear extension of $h_0$ to the
  whole $M_{\QQ}$, for $m\in\sigma^{\vee}_M$ we obtain
$$A_m=t^{-\lfloor h(m)\rfloor} \KK[t]=\bigoplus_{r\geq-h(m)}\KK t^r\,.$$

Letting $\widehat{N}=N\times \ZZ$, $\widehat{M}=M\times \ZZ$, and
$\widehat{\sigma}$ be the cone in $\widehat{N}_{\QQ}$ spanned by
$(\sigma,0)$ and $(p,1)$, a vector $(m,r)\in \widehat{M}_{\QQ}$
belongs to the dual cone $\widehat{\sigma}^{\vee}$ if and only if
$m\in\sigma^{\vee}$ and $r\geq-h(m)$. By identifying $\chi^{(0,1)}$
with $t$ we obtain
$$A=\bigoplus_{(m,r)\in\widehat{\sigma}^\vee_{\widehat{M}}}\KK t^r\chi^m=\bigoplus_{(m,r)\in\widehat{\sigma}^\vee_{\widehat{M}}}\KK\chi^{(m,r)} =\KK[\widehat{\sigma}_{\widehat{M}}]\,.$$
Hence $A$ is an affine semigroup algebra and so, we can apply the
results of the previous section.

Since $A_0$ is spanned as affine semigroup algebra by the character
$\chi^{(0,1)}$, the only codimension 1 face of $\widehat{\sigma}^\vee$
not containing the lattice vectors $(0,1)$
is $$\tau=\{(m,r)\in\widehat{M}_\QQ\mid m\in\sigma^\vee,\ r=-h(m)\}\,.$$
This is the face of $\widehat{\sigma}^\vee$ dual to the extremal ray
$\rho$ spanned by $(p,1)$ in $\widehat{N}_\QQ$.

In the notation of Lemma \ref{convex-core}, picking $e'\in S_\rho$ and
$\lambda\in\KK^*$ we let $\partial=\lambda\partial_{\rho,e'}$ be the
homogeneous LND with respect to the $\widehat{M}$-grading described in
Lemma \ref{toric-LND}.  Since $(0,1)\notin \tau$, $\partial$ is of
horizontal type with respect to the $M$-grading on $A$. Let $\deg_M$
stand for the corresponding degree function.

For any $e'=(e,s)\in M\times\ZZ$ we have $\deg_M\partial=e$ and
$\ker\partial=\KK[\tau_{\widehat{M}}]$. Therefore, in the notation of
Lemma \ref{kernel-linear}, $\omega^\vee=\sigma^\vee$ and $L=\{m\in
M\mid h(m)\in\ZZ\}$.

To be more concrete, we let $d>0$ be the smallest integer such that
$d\cdot p\in N$. Then $d\cdot h$ is an integer valued function on
$\sigma^\vee_M$. Letting $m_1\in M$ be a lattice vector such that
$\{h(m_1)\}=\{\frac{1}{d}\}$, by a routine calculation we obtain
\begin{align} \label{the-Sr} S_{\rho}=\left\{(e,s)\in \widehat{M}\mid
    e\in L-m_1,\ s=-h(e)-\tfrac{1}{d}\right\}\cap\sigma_\rho^{\vee}\,,
\end{align}
and
\begin{align} \label{the-derivation}
  \partial(\chi^m\cdot t^r)=\lambda\left(r+h(m)\right)\cdot\chi^{m+e}\cdot
t^{r-h(e)-\nicefrac{1}{d}},\quad \forall\ (m,r)\in\widehat{M}
\end{align}
where $\sigma_\rho\subseteq \widehat{N}_{\QQ}$ is as defined in Lemma
\ref{convex-core}, $\lambda\in\KK^*$, and $\partial_t$ is the partial
derivative with respect to $t$. Moreover, in this case
$\sigma_\rho=\sigma\times\{0\}$ and so
$$
S_{\rho}=\left\{(e,s)\in\widehat{M}\mid e\in \sigma^\vee\cap(L-m_1),\
  s=-h(e)-\tfrac{1}{d}\right\}\,.
$$
\end{example}

\begin{example} \label{ell-exa} Let $C=\PP^1$, $p\in N_{\QQ}$. Let
  $\Delta_\infty$ be a $\sigma$-tailed polyhedron (see Definition
  \ref{stp} (i)), and let $\DD=(p+\sigma)\cdot [0]+\Delta_\infty\cdot
  [\infty]$. Under these assumptions $h_0:\sigma^\vee\rightarrow\QQ,\,
  m\mapsto\langle m,p\rangle$ is linear and $h_z=0\ \forall
  z\in\KK^*$. We let as before $h:M_\QQ\rightarrow \QQ$ denote the
  linear extension of $h_0$ to the whole $M_{\QQ}$. We also suppose
  that $p+\Delta_\infty\subsetneq\sigma$ and so the sum
  $h_0+h_{\infty}\geq 0$ is not identically 0. Under these assumptions
  the $\sigma$-polyhedral divisor $\DD$ is proper in the sense of
  Definition \ref{ppd}. Letting $A=A[C,\DD]$, for any
  $m\in\sigma_M^\vee$ we have
$$A_m=\bigoplus_{-h_0(m)\leq r\leq h_{\infty}(m)}\KK t^r\,.$$

Let $\widehat{N}=N\times \ZZ$, $\widehat{M}=M\times \ZZ$, and let
$\widehat{\sigma}$ be the cone in $\widehat{N}_{\QQ}$ spanned by
$(\sigma,0)$, $(p,1)$ and $(\Delta_\infty,-1)$. A vector $(m,r)\in
\widehat{M}_{\QQ}$ belongs to the dual cone $\widehat{\sigma}^{\vee}$
if and only if $m\in\sigma^{\vee}$, $r\geq-h_0(m)$ and $r\leq
h_{\infty}(m)$. Thus by identifying $\chi^{(0,1)}$ with $t$ we obtain:
$$A=\bigoplus_{(m,r)\in\widehat{\sigma}^\vee_{\widehat{M}}}\KK t^r\chi^m=\bigoplus_{(m,r)\in\widehat{\sigma}^\vee_{\widehat{M}}}\KK\chi^{(m,r)} =\KK[\widehat{\sigma}_{\widehat{M}}]\,.$$
Hence $A$ is again an affine semigroup algebra, and so the results in
the previous section can be applied.

We let as before $\rho\subseteq\widehat{\sigma}$ be the extremal ray
spanned by $(p,1)$. The codimension 1 face dual to $\rho$ is
$$\tau=\{(m,r)\in\widehat{M}_\QQ\mid m\in\sigma^\vee,\ r=-h(m)\}\,.$$

In the notation of Lemma \ref{convex-core}, picking $e'\in S_\rho$ and
$\lambda\in\KK^*$ we let $\partial=\lambda\partial_{\rho,e'}$ be the
homogeneous LND with respect to the $\widehat{M}$-grading described in
Lemma \ref{toric-LND}.  Again $\partial$ is of horizontal type with
respect to the $M$-grading on $A$.

Furthermore, for any $e'=(e,s)\in M\times\ZZ$ we have
$\deg_M\partial=e$ and
$\ker\partial=\KK[\tau_{\widehat{M}}]$. Therefore, in the notation of
Lemma \ref{kernel-linear}, $\omega^\vee=\sigma^\vee$ and $L=\{m\in
M\mid h(m)\in\ZZ\}$.

To be more concrete, we let $d$ and $m_1$ be as in the previous
example. By a routine calculation we obtain that $S_\rho$ is as in
\eqref{the-Sr} and $\partial$ is as in \eqref{the-derivation}.
\end{example}

\begin{remark} \label{the-derivation-prop}
  \begin{enumerate}[(1)]
  \item In both examples, the homogeneous LND $\partial$ extends to a
    derivation on $K_0[M]$ given by \eqref{the-derivation}.
  \item With the same formula \eqref{the-derivation}, $\partial$
    extends to a homogeneous LND on $$A_M:=\bigoplus_{m\in
      M}t^{-\lfloor h(m)\rfloor} \KK[t]\chi^m,\quad\mbox{where}\quad
    A\subseteq A_M \subseteq K_0[M]\,.$$
  \item In particular, if $p=0$, then $\rho$ is the extremal ray
    spanned by $(0,1)$, $d=1$, and $L=M$. Furthermore, we can
    choose $m_1=0$ so that
    $S_{\rho}=(M\times\{-1\})\cap\sigma_1^\vee$, and the homogeneous
    LND $\partial$ of horizontal type on $A$ is given by
    $\partial=\lambda\chi^e\partial_t$, where $(e,-1)\in S_\rho$.
  \end{enumerate}
\end{remark}

We return now to the general case. We recall that
$$A=A[C,\DD],\quad \mbox{where}\quad \DD=\sum_{z\in C}\Delta_z\cdot z$$ is a proper $\sigma$-polyhedral divisor on $C=\AF^1$ or $C=\PP^1$, $h_z$ is the support function of $\Delta_z$, and $\partial$ is a homogeneous LND of horizontal type on $A$. 

In the next lemma we show that the subalgebra $A_\omega$ of $A$
generated by the homogeneous elements whose degrees are contained in
$\omega^\vee$, is as in the previous examples.

\begin{lemma} \label{hor-ome} With the notation of Lemma
  \ref{kernel-linear}, we let
  $A_{\omega}=\bigoplus_{m\in\omega^\vee_M}A_m\chi^m$. Then
  $A_{\omega}\simeq A[C,\DD_\omega]$ as $M$-graded algebras, where
  \begin{enumerate}[(i)]
  \item $\DD_\omega=(p+\omega)\cdot [0]$ for some $p\in N_{\QQ}$, in
    the case where $C=\AF^1$, and
  \item $\DD_\omega=(p+\omega)\cdot [0]+\Delta_\infty\cdot [\infty]$
    for some $p\in N_{\QQ}$ and some $\Delta_\infty\in
    \pol_{\sigma}(N_{\QQ})$ with $p+\Delta_\infty\subsetneq\sigma$, in
    the case where $C=\PP^1$.
  \end{enumerate}
\end{lemma}

\begin{proof}
  By Lemma \ref{kernel-linear} (3), the support functions $h_z$
  restricted to $\omega^\vee$ are linear for all $z\in \AF^1$ in the
  non-elliptic case and for all $z\in \PP^1\setminus\{z_\infty\}$ in
  the elliptic case. In the non-elliptic case this shows that
  $\DD_\omega=\sum_{z\in C} (p_z+\omega)\cdot z$, where $p_z\in
  N_{\QQ}$. In the elliptic case, we may suppose that
  $z_\infty=\infty$ and so $\DD_\omega=\sum_{z\in \AF^1}
  (p_z+\omega)\cdot z+\Delta_\infty\cdot[\infty]$, where
  $\Delta_\infty\in \pol_{\sigma}(N_{\QQ})$ and $p_z\in N_{\QQ}$
  $\forall z\in\AF^1$.

  By Lemma \ref{LND} (vi), without loss of generality we may assume
  that $\deg\partial\in\omega_M^\vee$. Letting $e=\deg\partial$ we
  consider the 2-dimensional finitely generated normal $\ZZ_{\geq
    0}$-graded domain
  $$B_e=\bigoplus_{r\in\ZZ_{\geq0}}A_{re}\chi^{re}\,.$$ 

  If $C$ is affine then $(B_e,\partial|_{B_e})$ is a parabolic pair in
  the sense of Definition 3.1 in \cite{FlZa2}. Now Corollary 3.19 in
  \emph{loc. cit.} shows that, for any $r\in\ZZ_{\geq0}$, the
  fractional part $\{\DD_\omega(re)\}$ is supported in at most one
  point\footnote{The classification results in \cite{FlZa2} are stated
    for surfaces over the field $\CC$ but they are valid over any
    algebraically closed field of characteristic 0 with the same
    proofs.}.  While for $C$ projective, $(B_e,\partial|_{B_e})$ is an
  elliptic pair in the sense of \emph{loc. cit.} Then Theorem 3.3 in
  \emph{loc. cit.} shows that $B_e$ is an affine semigroup
  algebra. According to Example 5.1 in \cite{Ti2}, for any
  $r\in\ZZ_{\geq0}$, the fractional part $\{\DD_\omega(re)\}$ is
  supported in at most two point.

Given $m\in L$, the derivation $\varphi_{m}\chi^m\partial$ on $A$ with
$\varphi_{m}$ as in Lemma \ref{kernel-linear} (1) is again locally
nilpotent. Applying the previous analysis to this LND shows that, for
any $r\in\ZZ_{\geq0}$, the fractional part
$\{\DD_\omega(r\cdot(e+m))\}$ is supported in at most one point in the
non-elliptic case and in at most two points in the elliptic case. By
Lemma \ref{kernel-linear} (4) $L$ and $e$ span $M$. So the functions
$h_z|_{\omega^\vee}$ are integral except for at most one value of $z$
in the non-elliptic case and at most two values of $z$ in the elliptic
case. Furthermore, in the elliptic case one of the two values of
$z\in\PP^1$ such that $h_z$ is not integral corresponds to $z=\infty$.

Without loss of generality, in both cases we may suppose that $z=0$ is
an exceptional value in $\AF^1$, provided there is one. In particular
$p_z\in N$ is a lattice vector for any $z\in\KK^*$. Since any integral
divisor on $\AF^1$ and any integral divisor of degree 0 on $\PP^1$ are
principal, Theorem \ref{AH-description} shows that $\DD_\omega$ can
always be chosen so that $p_z=0$ $\forall z\in\KK^*$. Now the result
follows.
\end{proof}

\begin{remark} \label{Lgood}
  \begin{enumerate}[(1)]
  \item By Examples \ref{hor-fund-exa} and \ref{ell-exa}, the previous
    lemma shows that $A_\omega$ is an affine semigroup algebra, or
    equivalently, $\spec\, A_\omega$ is a toric variety. Hence,
    $\spec\, A_\omega$ is a toric variety containing $X=\spec A$ as an
    open subset.
  \item In the notation of Lemma \ref{hor-ome}, let $h(m)=\langle
    m,p\rangle$. By virtue of Lemma \ref{kernel-linear} (1) and (2),
    $L=\{m\in M\mid h(m)\in\ZZ\}$.
  \end{enumerate}
\end{remark}

\begin{remark} \label{ker-tor} Whatever is an isomorphism $A\simeq
  A[C,\DD]$, the proof of the previous lemma implies the following.
  \begin{enumerate}[(1)]
  \item If $C=\AF^1$ then all $h_z|_{\omega^\vee}$ are linear and all
    but possibly one of them are integral.
  \item If $C=\PP^1$ then all but possibly one of $h_z|_{\omega^\vee}$
    are linear and all but possibly two of them are integral.
  \item By virtue of Theorem \ref{AH-description}, we may suppose, in
    both cases, that $h_z|_{\omega^\vee}=0$ $\forall z\in\KK^*$ and
    $h_0|_{\omega^\vee}$ is linear.
  \end{enumerate}
\end{remark}

The following lemma provides the main ingredient in our classification
of the homogeneous LNDs of horizontal type on $A=A[C,\DD]$.

\begin{lemma} \label{pres} Let $\DD$ be a proper $\sigma$-polyhedral
  divisor on $C=\AF^1$ or $C=\PP^1$. Let $\omega^{\vee}$ be a maximal
  cone in the quasifan $\Delta(\DD)$ or $\Delta(\DD|_{\AF^1})$,
  respectively, such that $h_z|_{\omega^\vee}=0\ \forall
  z\in\KK^*$. Let $\partial$ be the derivation of degree $e$ given by
  formula \eqref{the-derivation}. Then $\partial$ extends to a
  homogeneous LND on $A=A[C,\DD]$ if and only if, for every
  $m\in\sigma_M^\vee$ such that $m+e\in\sigma_M^\vee$ the following
  hold.
  \begin{enumerate}[(i)]
  \item If $h_z(m+e)\neq 0$, then $\lfloor h_z(m+e)\rfloor -\lfloor
    h_z(m)\rfloor \geq 1$ $\forall z\in\KK^*$.
  \item If $h_0(m+e)\neq h(m+e)$, then $\lfloor dh_0(m+e)\rfloor
    -\lfloor dh_0(m)\rfloor \geq 1+d h(e)$.
  \item If $C=\PP^1$, then $\lfloor d h_\infty(m+e)\rfloor -\lfloor d
    h_\infty(m)\rfloor \geq -1-d h(e)$.
  \end{enumerate}
  Here $h$ is the linear extension of $h_0|_{\omega^\vee}$ and $d>0$
  is the smallest integer such that $dh$ is integral.
\end{lemma}

\begin{proof}
  Similarly as in Example \ref{hor-fund-exa}, $h(m)=\langle
  m,p\rangle$ for some $p\in N_{\QQ}$. Since each $h_z$ is upper
  convex (see Definition \ref{stp} (ii)), $h_z(m)\leq 0$ for
  $z\in\KK^*$ and $h_0(m)\leq h(m)$. Letting $A_M=\bigoplus_{m\in
    M}\varphi_m \KK[t]\chi^m$, where $\varphi_m=t^{-\lfloor
    h(m)\rfloor}$ (see Remark \ref{the-derivation-prop}) we have
  $A\subseteq A_M$. By virtue of this remark $\partial$ extends to a
  homogeneous LND on $A_M$. We still denote by $\partial$ this
  extension. Thus $\partial$ extends to a homogeneous LND on $A$ if
  and only if $\partial$ stabilizes $A$.

  To show that $\partial$ stabilizes $A$, let us start with the
  simplest case where
  $h=0$.
  
  {\bfseries Case \pmb{$h=0$}.} In this case, Remark
  \ref{the-derivation-prop} (3) shows that $L=M$, $d=1$, and $r=-1$,
  and so $\partial=\lambda\chi^e\partial_t$. Furthermore, $h_z\leq 0\
  \forall z\in\AF^1$ and in the elliptic case $h_\infty\geq 0$. For
  any $m\in\sigma_M^\vee$ such that $m+e\in\sigma_M^\vee$, the
  conditions in the lemma can be reduced to
  \begin{enumerate}
  \item[\emph{(i$'$)}] If $h_z(m+e)\neq 0$, then $\lfloor
    h_z(m+e)\rfloor -\lfloor h_z(m)\rfloor \geq 1$ $\forall
    z\in\AF^1$.
  \item[\emph{(iii$'$)}] If $C=\PP^1$, then $\lfloor
    h_\infty(m+e)\rfloor -\lfloor h_\infty(m)\rfloor \geq -1$ $\forall
    m\in \sigma_M^\vee$.
  \end{enumerate}

  In this case
  $A_m=H^0\left(C,\OO(\lfloor\DD(m)\rfloor)\right)\subseteq \KK[t]$
  and $\partial$ stabilizes $A$ if and only if
$$f(t)\in A_m\Rightarrow f'(t)\in A_{m+e}\,, \forall m \in \sigma_M^\vee\,,$$
or equivalently
$$\divi\,f+\lfloor \DD(m)\rfloor\geq 0\Rightarrow \divi\, f'+\lfloor\DD(m+e)\rfloor\geq0\,,\forall m \in \sigma_M^\vee\,,$$
or else
\begin{align} \label{stab} \ord_z(f)+\lfloor h_z(m)\rfloor\geq 0
  \Rightarrow \ord_z(f')+\lfloor h_z(m+e)\rfloor\geq 0\,,\forall m \in
  \sigma_M^\vee\mbox{ and }\forall z\in C\,.
\end{align}

Next we show that (i$'$) and (iii$'$) hold if and only if \eqref{stab}
holds.

Let $z\in\AF^1$ and let $m\in\sigma_M^\vee$ such that
$m+e\in\sigma_M^\vee$.  If $h_z(m+e)=0$ the condition \eqref{stab}
holds since $f\in\KK[t]$.

Assume $h_z(m+e)\neq 0$. Since $h_z\leq 0$ is upper convex, if
$h_z(m)=0$ then $h_z(m+re)\neq 0$ $\forall r>1$ contradicting the fact
that $\partial$ is an LND. Hence we may assume that $h_z(m)\neq 0$ so
that $f\in(t-z)\KK[t]$. In this setting $\ord_z(f')=\ord_z(f)-1$ and
so
\begin{align} \label{stab2} \ord_z(f')+\lfloor
  h_z(m+e)\rfloor=\ord_z(f)+\lfloor h_z(m)\rfloor+\left( \lfloor
    h_z(m+e)\rfloor-\lfloor h_z(m)\rfloor-1\right)\,.
\end{align}
Therefore (i$'$) implies \eqref{stab}.

To show the converse, let us suppose that \eqref{stab} holds. Assuming
that $C$ is affine, for every $m\in\sigma_M^\vee$ we consider
$\varphi_m$ as in Lemma \ref{kernel-linear}. Since by this lemma
$\ord_z(\varphi_m)+\lfloor h_z(m)\rfloor=0$, applying \eqref{stab} and
\eqref{stab2} to $\varphi_m$ we obtain
$$\ord_z(\varphi_m)+\lfloor h_z(m)\rfloor+\left( \lfloor h_z(m+e)\rfloor-\lfloor h_z(m)\rfloor-1\right)=\lfloor h_z(m+e)\rfloor-\lfloor h_z(m)\rfloor-1\geq 0\,,$$
proving (i$'$) when $C$ is affine. If $C$ is projective, then for any
$z\in \AF^1$ and any $m\in\sigma_M^\vee$ we can still find
$\varphi_{m,z}\in A_{m}$ such that $\ord_z(\varphi_{m,z})+\lfloor
h_z(m)\rfloor=0$. Thus again the previous argument applies.

In the elliptic case, we let $z=\infty$ and we fix
$m\in\sigma_M^\vee$. If $f$ is constant, then \eqref{stab} holds
because $h_\infty(m)\geq 0$. Otherwise
$\ord_\infty(f')=\ord_\infty(f)+1$ and so
\begin{align}\label{stab3}
  \ord_\infty(f')+\lfloor h_\infty(m+e)\rfloor =\ord_\infty(f)+\lfloor
  h_\infty(m)\rfloor+\left( \lfloor h_\infty(m+e)\rfloor-\lfloor
    h_\infty(m)\rfloor+1\right)\,.
\end{align}
Therefore (iii$'$) implies \eqref{stab}.

To show the converse, we let as before $\varphi_{m,\infty}\in A_{m}$
be such that $\ord_\infty(\varphi_{m,\infty})+\lfloor
h_\infty(m)\rfloor=0$. Applying \eqref{stab} and \eqref{stab3} to
$\varphi_{m,\infty}$ we obtain
\begin{align*}
  \ord_\infty(\varphi_{m,\infty})+\lfloor h_\infty(m)\rfloor+\left(
    \lfloor h_\infty(m+e)\rfloor-\lfloor h_\infty(m)\rfloor+1\right)
  =\lfloor h_\infty(m+e)\rfloor-\lfloor h_{\infty}(m)\rfloor+1\geq
  0\,,
\end{align*}
proving (iii$'$).

Next we assume that $h$ is integral.

{\bfseries Case \pmb{$h$} integral.} In this case we still have
$d=1$. We recall that $h(m)=\langle m,p\rangle$. Letting
$\DD'=\DD-(p+\sigma)\cdot[0]$ if $C$ is affine and
$\DD'=\DD-(p+\sigma)\cdot[0]+(p+\sigma)\cdot[\infty]$ if $C$ is
projective, by Theorem \ref{AH-description} (iii) $A\simeq
A[C,\DD']$. In this setting $A[C,\DD']$ is as in the previous case
with $h'_0=h_0-h$, $h'_\infty=h_\infty+h$ and $h'_z=h_z$ $\forall z\in
\KK^*$.

This consideration shows that $\partial$ stabilizes $A$ if and only if
(i$'$) and (iii$'$) hold for $h'_z(m)\ \forall z\in C$. For any
$z\in\KK^*$, (i$'$) is equivalent to (i) in the lemma. Since
$$\lfloor h'_0(m+e)\rfloor -\lfloor h'_0(m)\rfloor-1=\lfloor h_0(m+e)\rfloor -\lfloor h_0(m)\rfloor-1-h(e)\,,$$
condition (i$'$) for $z=0$ is equivalent to (ii).

Similarly, if $C$ is projective
$$\lfloor h'_\infty(m+e)\rfloor -\lfloor h'_\infty(m)\rfloor+1=\lfloor h_\infty(m+e)\rfloor -\lfloor h_\infty(m)\rfloor+1+h(e)\,,$$
and so (iii$'$) is equivalent to (iii).

Now we turn to the general case.

{\bfseries General case.}  We may assume that $h$ is not integral
i.e., $d>1$. We consider the normalization $A'$ of
$A[\sqrt[d]{\varphi_{de}}\chi^{e}]$, where $\varphi_{de}:=t^{- h(de)}$
so that $A\subseteq A'$ is a cyclic extension.  With the notation of
Lemma \ref{cyclic-divisor} we have $A'=A[C',\DD']$ and
$K'_0=K_0[\sqrt[d]{\varphi_{de}}]$.

By the minimality of $d$ we deduce that $\gcd(h(de),d)=1$ and so
$\sqrt[d]{\varphi_{de}}=t^{a+b/d}$, where $\gcd(b,d)=1$. So
$K_0'=\KK(s)$, where $s^d=t$. Thus $C'\simeq\AF^1$ if $A$ is
non-elliptic and $C'\simeq\PP^1$ if $A$ is elliptic. Let
$p:C'\rightarrow C,\ z'\mapsto z'^d=z$ be the projection induced by
the morphism $K_0\hookrightarrow K'_0,\ t\mapsto t=s^d$. By Lemma
\ref{cyclic-divisor} we have
$$
\DD'=d\cdot\Delta_0\cdot [0]+\sum_{z'\in \KK^*}\Delta_{z}\cdot z'\
\text{if}\ C=\AF^1\,,$$ and
$$
\DD'=d\cdot\Delta_0\cdot [0]+d\cdot\Delta_\infty\cdot
[\infty]+\sum_{z'\in \KK^*}\Delta_{z}\cdot z'\ \text{if}\ C=\PP^1\,.
$$
So $h'_0=d h_0$, $h'_\infty=d h_\infty$ and $h'_{z'}=h_{z}$. Moreover
$h'_0|_{\omega^\vee}$ is integral and $A'$ is as in the previous case.

Recall that $A_M=\bigoplus_{m\in M}\varphi_m \KK[t]\chi^m$, where
$\varphi_m=t^{-\lfloor h(m)\rfloor}$. We define further
$$A_M'=\bigoplus_{m\in M}\varphi'_m \KK[s]\chi^m,\qquad \mbox{where}\qquad \varphi'_m=-s^{dh(m)}\,.$$
Since $A_M\subseteq A_M'$ is a cyclic extension, by Lemma
\ref{cyclic-LND}, $\partial:A_M\rightarrow A_M$ extends to a
homogeneous LND $\partial':A'_M\rightarrow A'_M$.

We claim that $\partial$ stabilizes $A$ if and only if $\partial'$
stabilizes $A'$. In fact the ``only if'' direction is a consequence of
Lemma \ref{cyclic-LND}. If $\partial'$ stabilizes $A'$ then
$\partial'(A)=\partial(A)\subseteq A_M\cap A'=A$, proving the claim.

We let $h'$ be the linear extension of $h'_0|_{\omega^\vee}$. Clearly
$h'=dh$. The previous case shows that $\partial'$ stabilizes $A'$ if
and only if, for any $m\in\sigma_M^\vee$ such that
$m+e\in\sigma_M^\vee$, the following conditions hold.
\begin{enumerate}
\item[\emph{(i$''$)}] If $h'_{z'}(m+e)\neq 0$, then $\lfloor
  h'_{z'}(m+e)\rfloor -\lfloor h'_{z'}(m)\rfloor \geq 1$ $\forall
  z'\in\KK^*$.
\item[\emph{(ii$''$)}] If $h'_0(m+e)\neq h'(m+e)$, then $\lfloor
  h'_0(m+e)\rfloor -\lfloor h'_0(m)\rfloor \geq 1+h'(e)$.
\item[\emph{(iii$''$)}] If $C=\PP^1$, then $\lfloor
  h'_\infty(m+e)\rfloor -\lfloor h'_\infty(m)\rfloor \geq -1-h'(e)$.
\end{enumerate}

Replacing in (i$''$)-(iii$''$) $h'$ by $dh$, $h'_0$ by $dh_0$,
$h'_\infty$ by $dh_\infty$, and $h'_{z'}$ by $h_{z}$ for $z\in\KK^*$,
shows that $\partial$ stabilizes $A$ if and only if (i)-(iii) of the
lemma hold. Now the proof is completed.
\end{proof}

\begin{remark} \label{iiisup} In the elliptic case, if
  $e\in\omega_M^\vee$, then (iii) in Lemma \ref{pres} holds. In fact
  \begin{align*}
    \lfloor d h_\infty(m+e)\rfloor -\lfloor d h_\infty(m)\rfloor&\geq dh_\infty(m+e)-1-dh_\infty(m)\\
    &\geq dh_\infty(e)-1\geq -d h(e)-1\,.
  \end{align*}
\end{remark}

In the following theorem we describe all the homogeneous LND of
horizontal type on a $\TT$-variety of complexity one. It is our main
classification result which summarizes the previous ones.

\begin{theorem} \label{resultado} Let $\DD$ be a proper
  $\sigma$-polyhedral divisor on $C=\AF^1$ or $C=\PP^1$, and let
  $A=A[C,\DD]$. Let $\omega^\vee\subseteq M_{\QQ}$ be a polyhedral
  cone, and $e\in M$ be a lattice vector. Then there exists a
  homogeneous LND $\partial:A\rightarrow A$ of horizontal type with
  $\deg\partial=e$ and $\omega^\vee(\partial)=\omega^\vee$ if and only
  if the following conditions (i)-(v) hold.
  \begin{enumerate}
  \item[(i)] If $C=\AF^1$, then $\omega^\vee$ is a maximal cone in the
    quasifan $\Lambda(\DD)$, and there exists $z_0\in C$ such that
    $h_z|_{\omega^\vee}$ is integral $\forall z\in C\setminus\{z_0\}$.
  \item[(i$'$)] If $C=\PP^1$, then there exists $z_\infty\in\PP^1$
    such that (i) holds for $C_0:=\PP^1\setminus\{z_\infty\}$.
  \end{enumerate}
  Without loss of generality, we may suppose that $z_0=0$,
  $z_\infty=\infty$ in the elliptic case, and
  $h_z(m)|_{\omega^\vee}=0$ $\forall z\in\KK^*$. Let $h$ and $d$ be as
  in Lemma \ref{pres}, let $m_1$ be as in Example \ref{hor-fund-exa},
  and let $L$ be as in Remark \ref{Lgood} (2).
  \begin{enumerate}
  \item[(ii)] The lattice vector $(e,-\tfrac{1}{d}-h(e))$ belongs to
    $S_\rho$ as defined in \eqref{the-Sr}.
  \end{enumerate}
  For any $m\in\sigma_M^\vee$ such that $m+e\in\sigma_M^\vee$, the
  following hold.
  \begin{enumerate}
  \item[(iii)] If $h_z(m+e)\neq 0$, then $\lfloor h_z(m+e)\rfloor
    -\lfloor h_z(m)\rfloor \geq 1$ $\forall z\in\KK^*$.
  \item[(iv)] If $h_0(m+e)\neq h(m+e)$, then $\lfloor dh_0(m+e)\rfloor
    -\lfloor dh_0(m)\rfloor \geq 1+d h(e)$.
  \item[(v)] If $C=\PP^1$, then $\lfloor d h_\infty(m+e)\rfloor
    -\lfloor d h_\infty(m)\rfloor \geq -1-d h(e)$.
  \end{enumerate}
  Moreover,
$$\ker\partial=\bigoplus_{m\in\omega_L^\vee}\KK\varphi_m\chi^m\,,$$
where $\varphi_m\in A_m$ satisfy the relation
$$\divi(\varphi_m)+\DD(m)=0 \mbox{\ \ if\ \ } C=\AF^1 \quad \mbox{or} \quad \divi(\varphi_m)|_{C_0}+\DD(m)|_{C_0}=0 \mbox{\ \ if\ \ } C=\PP^1\,.$$
\end{theorem}

\begin{proof}
  Let $\partial$ be a homogeneous LND of horizontal type on $A$ with
  $\deg\partial=e$ and $\omega^\vee(\partial)=\omega^\vee$. Lemma
  \ref{kernel-linear} (3) and Remark \ref{ker-tor} show that (i) and
  (i$'$) hold. Lemma \ref{hor-ome} and Examples \ref{hor-fund-exa} and
  \ref{ell-exa} shows that (ii) holds. To conclude, Lemma \ref{pres}
  shows that (iii)-(v) hold.

  To show the converse, assume that (i), (i$'$) and (ii)-(v) are
  fulfilled. By Theorem \ref{AH-description}, (i) and (i$'$) imply
  that $A_\omega\simeq A[C,\DD_\omega]$ with $\DD_\omega$ as in Lemma
  \ref{hor-ome}. By Examples \ref{hor-fund-exa} and \ref{ell-exa} and
  Remark \ref{the-derivation-prop} (2), (ii) shows that there exists a
  homogeneous LND $\partial:A_M\rightarrow A_M$ with
  $\deg\partial=e$. By Lemma \ref{pres} and its proof, (iii)-(v) imply
  that $\partial$ restricts to a homogeneous LND on $A$. Finally, by
  Lemma \ref{kernel-linear} (3), (i) and (i$'$) imply that
  $\omega^\vee(\partial)=\omega^\vee$.

  Moreover, Lemma \ref{kernel-linear} (1) and (2) give the desired
  description of $\ker\partial$.
\end{proof}

\begin{remark}
  The maximal cones in the quasifan $\Lambda(\DD)$ are in one to one
  correspondence with the vertices of the $\sigma$-polyhedron
  $\deg\DD$.
\end{remark}

\begin{corollary} \label{rege} In the notation of Theorem
  \ref{resultado}, $A$ admits a homogeneous LND $\partial$ of
  horizontal type such that $\omega^\vee(\partial)=\omega^\vee$ if and
  only if (i) and (i$'$) in the theorem hold.
\end{corollary}

\begin{proof}
  The ``only if'' part follows directly form Theorem \ref{pres}.

  Assume that (i) and (i$'$) hold. By Theorem \ref{pres} and Examples
  \ref{hor-fund-exa} and \ref{ell-exa}, we only need to show that
  there exists $e\in M$ such that
  $\left(e,-\tfrac{1}{d}-h(e)\right)\in S_\rho$ and (iii)-(v) hold.

  Let $(e',r')\in S_\rho$ (by Remark \ref{toric-non-empty}, this set
  is non-empty). By this remark $e=e'+m\ \forall m\in\omega_L^\vee$ is
  such that $(e,r'-h(m))\in S_\rho$. In particular, we can assume that
  $e$ belongs to the relative interior of $\omega^\vee$. In this
  setting, Remark \ref{iiisup} shows that (v) holds.

  As in the proof of Lemma \ref{3.4}, for every $z\in \AF^1$, we let
  $\{\delta_{0,z},\cdots,\delta_{\ell_z,z}\}$ denote the set of all
  maximal cones in $\Lambda(h_z)$ and $g_{r,z},\
  r\in\{0,\cdots,\ell_z\}$ be the linear extension of
  $h_z|_{\delta_{r,z}}$ to $M_{\QQ}$. We assume further that
  $\omega^\vee\subseteq\delta_{0,z}\,\forall z\in\AF^1$.

  Since the functions $h_z$ are upper convex, the inequalities in
  (iii) and (iv) hold if they hold in every maximal cone on
  $\Lambda(h_z)$ except $\delta_{0,z}$ i.e.,
  \begin{enumerate}
  \item[(iii$'$)] $\lfloor g_{r,z}(m+e)\rfloor -\lfloor g_z(m)\rfloor
    \geq 1$ $\forall z\in\KK^*$, $\forall r\in \{1,\cdots,\ell_z\}$
    and $\forall m\in\delta_{r,z}\cap M$.
  \item[(iv$'$)] $\lfloor dg_{r,0}(m+e)\rfloor -\lfloor
    dg_{r,0}(m)\rfloor \geq 1+d h(e)$ $\forall r\in
    \{1,\cdots,\ell_0\}$ and $\forall m\in\delta_{r,0}\cap M$.
  \end{enumerate}

  These inequalities are fulfilled if
  \begin{align} \label{eq7} g_{r,z}(e)\geq 1\ \forall z\in\KK^*\mbox{
      and } \forall r\in\{1,\cdots,\ell_z\}, \qquad\mbox{and}\qquad
    g_{r,0}(e)\geq \tfrac{1}{d}+\lceil h(e)\rceil\ \forall
    r\in\{1,\cdots,\ell_0\}\,.
  \end{align}

  Since $e$ belongs to the relative interior of $\omega^\vee$, we have
  $g_{r,z}(e)>g_{0,z}(e)\ \forall z\in\AF^1$, $g_{0,0}(e)=h(e)$, and
  $g_{0,z}=0\ \forall z\in\KK^*$. By the linearity of the functions
  $g_{r,z}$ we can choose $e$ such that \eqref{eq7} holds, proving the
  corollary.
\end{proof}

\begin{corollary} \label{hor-eq} In the notation on Theorem
  \ref{resultado}, two homogeneous LND $\partial$ and $\partial'$ of
  horizontal type on $A$ are equivalent if and only if
  $\omega^\vee(\partial)=\omega^\vee(\partial')$ and, in the elliptic
  case, $z_\infty(\partial)=z_\infty(\partial')$.
\end{corollary}
\begin{proof}
  Indeed, the description of $\ker\partial$ given in Theorem
  \ref{resultado} depends only on $\omega^\vee$ in the non-elliptic
  case and on $\omega^\vee$ and $z_\infty\in C$ in the elliptic one.
\end{proof}

\begin{corollary} \label{hor-equiv} The number of pairwise
  non-equivalent homogeneous LNDs of horizontal type on $A=A[C,\DD]$
  is finite except in the case where $A$ is elliptic and there exists
  a maximal cone $\omega^\vee$ of $\Lambda(\DD)$ such that all but
  possibly one $h_z|_{\omega^\vee}$ are integral.
\end{corollary}
\begin{proof}
  Since $\Lambda(\DD)$ has only a finite number of maximal cones,
  Corollary \ref{hor-eq} gives the result in the case where $A$ is
  non-elliptic. Furthermore, in the elliptic case by this corollary
  there is an infinite number of pairwise non-equivalent LNDs on $A$
  if and only if in Theorem \ref{resultado} (i$'$) we can choose
  $z_\infty\in\PP^1$ arbitrarily. However the latter is indeed
  possible under the assumptions of the corollary.
\end{proof}

\begin{example}
  A combinatorial description of $\KK^{[2]}=\KK[x,y]$ with the grading
  induced by $\deg x=\deg y=1$ is given by the proper
  $\sigma$-polyhedral divisor $\DD=(1+\sigma)\cdot[0]$ on $\PP^1$,
  where $\sigma=\QQ_{\geq0}\subseteq N_{\QQ}\simeq\QQ$. By Corollary
  \ref{hor-equiv} there exist an infinite number of pairwise
  non-equivalent LNDs on $\KK^{[2]}$ homogeneous with respect to the
  given grading. Indeed, the derivations on the family
$$\partial_{\lambda}=\lambda\frac{\partial}{\partial x}+(1-\lambda)\frac{\partial}{\partial y}$$
are homogeneous and pairwise non-equivalent for different values of
$\lambda$.

In contrast, a combinatorial description of $\KK^{[2]}$ with the
grading induced by $\deg x=-\deg y=1$ is given by the proper
$\sigma$-polyhedral divisor $\DD=[0,1]\cdot[0]$ on $\AF^1$. By
Corollary \ref{hor-equiv} there exist a finite number of pairwise
non-equivalent LNDs homogeneous with respect to this grading. Indeed,
by Corollary \ref{rege} the only such LNDs are the partial
derivatives.
\end{example}

\begin{remark}
  Let $A$ be a normal finitely generated effectively $M$-graded
  algebra, such that the complexity of the corresponding $\TT$-action
  on $\spec\, A$ is 0 or 1. In Corollaries \ref{fg-toric} and
  \ref{fgraro}, and Proposition \ref{fg-fib} we have shown that the
  kernel of a homogeneous LND on $A$ is finitely generated.

  On the other hand, there are examples of homogeneous LNDs on $\AF^r$
  for $r\geq 5$, whose kernel is not finitely generated, see
  \cite{Ro}, \cite{Fr1} and \cite{DaFr}. For instance, Daigle and
  Freudenburg showed in \cite{DaFr} that $\ker\partial$ is not
  finitely generated for the
  LND $$\partial=x_1^3\frac{\partial}{\partial
    x_2}+x_2\frac{\partial}{\partial x_3}+x_3\frac{\partial}{\partial
    x_4}+x_1^2\frac{\partial}{\partial x_5}$$ on
  $\KK^{[5]}=\KK[x_1,\ldots,x_5]$. Furthermore it is easy to see that
  $\partial$ is homogeneous of degree $(0,-1)$ under the effective
  $\ZZ^2$-grading on $\KK^{[5]}$ given by
$$\deg x_1=(1,0),\quad\deg x_2=(3,1),\quad\deg x_3=(3,2),\quad\deg x_4=(3,3),\quad\deg x_5=(2,1)\,.$$
The corresponding $\TT$-action on $\AF^5$ is of complexity 3.
\end{remark}

In the following example we study the existence of homogeneous LNDs on
the $M$-graded algebra $A$ of Example \ref{ex-hyp}.
\begin{example} \label{ex-hyp-rev} Let the notation be as in Example
  \ref{ex-hyp}. Since $\sigma=\{0\}$, Lemma \ref{fiber-negative} shows
  that there is no homogeneous LND of fiber type on $A$. In contrast,
  let us show that there exist exactly 4 pairwise non-equivalent
  homogeneous LNDs on $A$.

  Indeed, since $h_0$ is the only support function which is
  non-integral Corollaries \ref{rege} and \ref{hor-eq} show that there
  are four non-equivalent homogeneous LNDs of horizontal type on $A$
  corresponding to the four maximal cones in $\Lambda(\DD)$,
$$
\begin{array}{cc}
  \delta_1=\cone((1,0),(-4,1)), & \delta_2=\cone((-4,1),(-1,0)), \\
  \delta_3=\cone((-1,0),(8,-1)), & \delta_4=\cone((8,-1),(1,0))\,.
\end{array}
$$
For the cones $\delta_1$ and $\delta_2$ the hypothesis of Lemma
\ref{pres} are fulfilled i.e., $h_z|_{\delta_i}=0$ $\forall z\in\KK^*$
for $i=1,2$. Moreover, $e_1=(-3,1)$ and $e_2=(-8,1)$ satisfy
conditions (i)-(iii) in this lemma for $\delta_1$ and $\delta_2$,
respectively.

We let $\partial_1$ and $\partial_2$ be the respective LNDs defined in
\eqref{the-derivation}. Letting $m=(m_1,m_2)\in M$, by a routine
calculation we obtain
$$\partial_1\left(\chi^m t^r)=(r-\tfrac{1}{4}m_1-m_2\right)\cdot \chi^{m+e_1}t^r,
\qquad\mbox{and}\qquad
\partial_2(\chi^m t^r)=r\cdot\chi^{m+e_2}t^r\,.
$$

Furthermore, under the isomorphism \eqref{iso-ex1} in Example
\ref{ex-hyp}, $\partial_1$ and $\partial_2$ can be extended to
$\KK^{[4]}=\KK[x_1,x_2,x_3,x_4]$ as LNDs
$$\partial_1=-\frac{1}{4}x_3\frac{\partial}{\partial x_2}+x_1^2x_2^3\frac{\partial}{\partial x_4} \qquad\mbox{and}\qquad \partial_2=x_3\frac{\partial}{\partial x_1}-(2x_1x_2^4+1)\frac{\partial}{\partial x_4}\,.$$

To obtain the derivations corresponding to $\delta_3$ and $\delta_4$
we let $C'=\spec\,\KK[s]$, $\Delta'_1=\{0\}\times[-1,0]$, and
$\DD'=\Delta_0\cdot[0]+\Delta'_1\cdot[1]$. Theorem
\ref{AH-description} (3) shows that $A\simeq A[C',\DD']$. Under this
new combinatorial description we have
$$u_1=-s\chi^{(4,0)},\quad u_2=\chi^{(-1,0)},\quad u_3=(1-s)\chi^{(-4,1)}, \quad \mbox{and} \quad u_4=s\chi^{(8,-1)}\,.$$
Now the assumptions of Lemma \ref{pres} are satisfied for $\delta_3$
and $\delta_4$. Moreover, $e_3=(4,-1)$ and $e_4=(9,-1)$ satisfy
conditions (i)-(iii) in this lemma for $\delta_3$ and $\delta_4$,
respectively.

We let $\partial_3$ and $\partial_4$ be the respective LNDs defined by
\eqref{the-derivation}. By a simple computation we obtain
$$\partial_3\left(\chi^m s^r)=(r+m_2\right)\cdot \chi^{m+e_3}s^r,
\qquad\mbox{and}\qquad
\partial_4(\chi^m
s^r)=\left(r-\tfrac{1}{4}m_1-m_2\right)\cdot\chi^{m+e_4}s^{r+1}\,.
$$

Furthermore, under the isomorphism \eqref{iso-ex1} $\partial_3$ and
$\partial_4$ are induced by the LNDs
$$\partial_3=-x_4\frac{\partial}{\partial x_1}+(2x_1x_2^4+1)\frac{\partial}{\partial x_3} \qquad\mbox{and}\qquad \partial_4=\frac{1}{4}x_4\frac{\partial}{\partial x_2}-x_1^2x_2^3\frac{\partial}{\partial x_3}$$
on $\KK^{[4]}$.
\end{example}

\subsection{The surface case}
\label{sur-cas}

A description of $\CC^*$-surfaces was given in \cite{FlZa1} in terms
of the DPD (Dolgachev-Pinkham-Demazure) presentation. In \cite{FlZa2}
this description was applied to classify the homogeneous LNDs on
normal affine $\CC^*$-surfaces (of both horizontal and fiber
type). Here we relate both descriptions. Besides, we stress the
difference that appears in higher dimensions.

In the case of dimension 2 the lattice $N$ has rank 1, which makes
things quite explicit (cf. e.g., \cite{Su}).

We treat the elliptic case first. In this case $\sigma$ is of full
dimension, and so we can assume that $\sigma=\QQ_{\geq0}\subseteq
N_{\QQ}=\QQ$. Let $A=A[C,\DD]$, where $\DD$ is a proper
$\sigma$-polyhedral divisor on a smooth projective curve $C$. In this
setting, $\DD$ is uniquely determined by the $\QQ$-divisor $\DD(1)$ on
$C$. Here $(C,\DD(1))$ coincides with the DPD presentation data. Since
the only extremal ray of $\sigma$ is $\sigma$ itself and $\deg\DD$ is
$\sigma$-tailed (see Definition \ref{stp}), by Corollary \ref{fiex}
there is no homogeneous LND of fiber type on $A$.

Furthermore, if there is a homogeneous LND $\partial$ of horizontal
type on $A$, then $\omega^\vee(\partial)=\sigma^\vee$, and so by
Remark \ref{Lgood} (1) $A=A_{\omega}$ is an affine semigroup algebra
i.e., $\spec\, A$ is an affine toric surface. This corresponds to
Theorem 3.3 in \emph{loc. cit.}

Next we consider a non-elliptic algebra $A$ so that $C$ is an affine
curve. In \emph{loc.cit.} this case is further divided into two
subcases, the parabolic one which corresponds to $\sigma=\QQ_{\geq0}$,
and the hyperbolic one which corresponds to $\sigma=\{0\}$.

In the parabolic case, the DPD presentation data is the same as in the
elliptic one. In this case there is again just one extremal ray
$\rho=\sigma$ and $S_\rho=\{-1\}$. Moreover, since the support
functions $h_z$ are positively homogeneous on
$\sigma^\vee=\QQ_{\geq0}$, they are linear and so $D_{-1}=\DD(1)$ (see
Lemma \ref{3.4}). By Theorem \ref{fiber-clas} the homogeneous LNDs of
fiber type on $A$ are in one to one correspondence with the rational
functions
$$\varphi\in H^0(C,\OO_C(\lfloor-\DD(1)\rfloor))\,.$$
This corresponds to Theorem 3.12 in \emph{loc. cit.}

If a graded parabolic 2-dimensional algebra $A$ admits a homogeneous
LND of horizontal type, then $\spec\, A$ is a toric variety by the
same argument as in the elliptic case. This yields Theorem 3.16 and
Corollary 3.19 in \emph{loc. cit.}

In the hyperbolic case $\DD$ is uniquely determined by the pair of
$\QQ$-divisors $\left(\DD(1),\DD(-1)\right)$ which correspond to the
pair $(D_+,D_-)$ in the DPD presentation data. According to our
Definition \ref{stp} (ii), this pair satisfies $\DD(1)+\DD(-1)\leq
0$. In this case, by Lemma \ref{fiber-negative} there is no
homogeneous LND of fiber type on $A$ since $\sigma=\{0\}$. This
corresponds to Lemma 3.20 in \emph{loc. cit.}

The homogeneous LNDs of horizontal type are classified in Theorem
\ref{resultado} above. Specializing this classification to dimension 2
gives Theorem 3.22 in \emph{loc. cit.} More precisely, conditions (i)
and (ii) of \ref{resultado} lead to (i) of Theorem 3.22 in
\emph{loc. cit.} while (iii) and (iv) in \ref{resultado} lead to (ii)
in Theorem 3.22 in \emph{loc. cit.}

In contrast, in dimension 3 a new phenomena appear. For instance,
there exist non-toric threefolds with an elliptic $\TT$-action and a
homogeneous LND of horizontal or fiber type, see subsection
\ref{non-rat} for an example of fiber type. With the notation as in
subsection \ref{non-rat}, considering $C=\PP^1$ and
$\DD=\tfrac{1}{2}\Delta\cdot[0]+\tfrac{1}{2}\Delta\cdot[1]+\Delta'\cdot[\infty]$,
where $\Delta'=\sigma\cap\{\langle(1,1),\cdot\rangle\geq
1\}\subseteq N_{\QQ}$ gives a non-toric example with 2 equivalence
classes of homogeneous LNDs fiber type and 4 equivalence classes of
homogeneous LNDs of horizontal type.

\section{Applications}
\label{app}

\subsection{The Makar-Limanov invariant}

Let $A$ be a finitely generated normal domain, and let $\LND(A)$ be
the set of all LNDs on $A$. The \emph{Makar-Limanov invariant} of $A$
is defined as $$\ML(A)=\bigcap_{\partial\in\LND(A)}\ker\partial\,.$$
Similarly, if $A$ is effectively $M$-graded we let $\LND_h(A)$ be the
set of all homogeneous LNDs on $A$, and we call
$$\ML_h(A)=\bigcap_{\partial\in\LND_h(A)}\ker\partial$$
the \emph{homogeneous Makar-Limanov invariant} of $A$. Clearly
$\ML(A)\subseteq\ML_h(A)$.

In the sequel we apply the results in Section 2 and 3 in order to
compute $\ML_h(A)$ in the case where the complexity of the
$\TT$-action on $\spec\,A$ is 0 or 1. We also give some partial
results for the usual invariant $\ML(A)$ in this particular case.

\begin{remark}
  Since two equivalent LNDs (see Definition \ref{LND-equiv}) have the
  same kernel, to compute $\ML(A)$ or $\ML_h(A)$ it is sufficient to
  consider pairwise non-equivalent LNDs on $A$. The pairwise
  non-equivalent homogeneous LNDs on $A$ are classified in Corollary
  \ref{toric-neq} for complexity 0 case, and in Corollaries \ref{feq}
  and \ref{hor-eq} for complexity 1 case.
\end{remark}

We treat first the case of complexity 0 i.e., the case of affine toric
varieties. Let $\sigma\subseteq N_{\QQ}$ be a pointed polyhedral cone.

\begin{proposition} \label{toric-ML} Let $A=\KK[\sigma_M^\vee]$ be an
  affine semigroup algebra so that $X=\spec\,A$ is a toric
  variety. Then
$$
\ML(A)=\ML_h(A)=\KK[\theta_M]\,,
$$
where $\theta\subseteq M_{\QQ}$ is the maximal subspace contained in
$\sigma^{\vee}$. In particular $\ML(A)=\KK$ if and only if $\sigma$ is
of complete dimension i.e., if and only if there is no torus factor in
$X$.
\end{proposition}

\begin{proof}
  By Corollary \ref {toric-neq} and Theorem \ref{toric-clasification},
  the pairwise non-equivalent homogeneous LNDs on $A$ are in one to
  one correspondence with the extremal rays of $\sigma$. For any
  extremal ray $\rho\subseteq\sigma$ and any $e\in S_\rho$ as in Lemma
  \ref{convex-core}, the kernel of the corresponding homogeneous LND
  is $\ker\partial_{\rho,e}=\KK[\tau_M]$, where
  $\tau\subseteq\sigma^{\vee}$ is the codimension 1 face dual to
  $\rho$.

  Since $\theta\subseteq \sigma^\vee$ is the intersection of all
  codimension 1 faces, we have $\ML_h(A)=\KK[\theta_M]$. Furthermore,
  the characters in $\KK[\theta_M]\subseteq A$ are invertible
  functions on $A$ and so, by Lemma \ref{LND} (iii),
  $\partial(\KK[\theta_M])=0\ \forall \partial\in\LND(A)$. Hence
  $\KK[\theta_M]\subseteq \ML(A)$, proving the lemma.
\end{proof}

For the rest of this section, we let $A=A[C,\DD]$, where $\DD$ is a
proper $\sigma$-polyhedral divisor on a smooth curve $C$. We also let
$\ML_{fib}(A)$ and $\ML_{hor}(A)$ be the intersection of the kernels
of all homogeneous LNDs of fiber type and of horizontal type,
respectively, so that
\begin{align} \label{ML-eq} \ML_h(A)=\ML_{fib}(A)\cap\ML_{hor}(A)\,.
\end{align}

We first compute $\ML_{fib}(A)$. If $A$ is non-elliptic (elliptic,
respectively) we let $\{\rho_i\}$ be the set of all extremal rays of
$\sigma^\vee$ (of all extremal rays of $\sigma^\vee$ such that
$\rho\cap\deg\DD=\emptyset$, respectively). In both cases we let
$\tau_i\subseteq M_{\QQ}$ denote the codimension 1 face dual to
$\rho_i$ and $\theta=\bigcap\tau_i$.

\begin{lemma} \label{ML-fib} With the notation as above,
$$\ML_{fib}(A)=\bigoplus_{m\in \theta_M}A_m\chi^m\,.$$
\end{lemma}
\begin{proof}
  By Corollary \ref{fiex}, for every extremal ray $\rho_i$ there is a
  homogeneous LND $\partial_i$ of fiber type with kernel
  $\ker\partial_i=\bigoplus_{m\in\tau_i\cap M}A_m\chi^m$. By Corollary
  \ref{feq} any homogeneous LND of fiber type on $A$ is equivalent to
  one of the $\partial_i$. Finally, taking the intersection
  $\bigcap_i\ker\partial_i$ gives the desired description of
  $\ML_{fib}(A)$.
\end{proof}

\begin{remark} \label{int-fib} If $A$ is non-elliptic, then
  $\theta\subseteq M_{\QQ}$ is the maximal subspace contained in
  $\sigma^\vee$, as in the toric case. In particular, if $A$ is
  parabolic then $\theta=\{0\}$ and $\ML_{fib}(A)=A_0$, and if $A$ is
  hyperbolic then $\theta=M_{\QQ}$ and $\ML_{fib}(A)=A$.
\end{remark}

If there is no LND of horizontal type on $A$, then $\ML_{hor}(A)=A$
and $\ML_h(A)=\ML_{fib}(A)$. In the sequel we assume that $A$ admits a
homogeneous LND of horizontal type.

If $A$ is non-elliptic, we let $\{\delta_i\}$ be the set of all cones
in $M_{\QQ}$ satisfying (i) in Theorem \ref{resultado}, and
$\delta=\bigcap_i\delta_i$. If $A$ is elliptic, we let
$\{\delta_{i,z}\}$ be the set of all cones in $M_{\QQ}$ satisfying
(i$'$) in Theorem \ref{resultado} with $z_\infty=z$,
$B=\{m\in\sigma^\vee\mid h_{\deg\DD}=0\}$, and
$\delta=\bigcap_{i,z}\delta_{i,z}\cap B$.

\begin{lemma} \label{ML-hor} With the notation as before, if
  $\partial$ is a homogeneous LND on $A$ of horizontal type, then
$$\ML_{hor}(A)=\bigoplus_{m\in \delta_L}\KK\varphi_m\chi^m\,,$$
where $L=L(\partial)$ and $\varphi_m\in A_m$ satisfy the relation
$\divi(\varphi_m)+\DD(m)=0$.
\end{lemma}

\begin{proof}
  We treat first the non-elliptic case. By Corollary \ref{rege} for
  every $\delta_i$ there is a homogeneous LND $\partial_i$ of
  horizontal type with kernel
$$\ker\partial_i=\bigoplus_{m\in\delta_i\cap L_i}\KK\varphi_m\chi^m\,,$$
where $L_i=L(\partial_i)$ and $\varphi_m\in A_m$ is such that
$\divi(\varphi_m)+\DD(m)=0$. By Corollary \ref{hor-eq}, any
homogeneous LND of horizontal type on $A$ is equivalent to one of the
$\partial_i$. Taking the intersection of all $\ker\partial_i$ gives
the lemma in this case.

Let further $A$ be elliptic, and let $\partial$ be a homogeneous LND
of horizontal type on $A$. Let $z_0,z_\infty\in\PP_1$, and
$\omega^\vee$ and $L$ be as in Theorem \ref{resultado} so that
$$\ker\partial=\bigoplus_{m\in\omega^\vee_L}\KK\varphi_m\chi^m\,,$$ 
where $\varphi_m\in A_m$ satisfies
$\divi(\varphi_m)|_{\PP^1\setminus\{z_\infty\}}+\DD(m)|_{\PP^1\setminus\{z_\infty\}}=0$.

By permuting the roles of $z_0$ and $z_\infty$ in Theorem
\ref{resultado} we obtain another LND $\partial'$ on $A$. The
description of $\ker\partial$ and $\ker\partial'$ shows that
$$\ker\partial\cap\ker\partial'=\bigoplus_{\omega^\vee_L\cap B}\KK\varphi\chi^m\,,$$
where $\varphi_m\in A_m$ is such that $\divi(\varphi_m)+\DD(m)=0$.

Now the lemma follows by an argument similar to that in the
non-elliptic case.
\end{proof}

\begin{theorem} \label{ML-cp1} In the notation of Lemmas \ref{ML-fib}
  and \ref{ML-hor}, if there is no homogeneous LND of horizontal type
  on $A$, then
$$\ML_{h}(A)=\bigoplus_{m\in \theta_M}A_m\chi^m\,.$$
If $\partial$ is a homogeneous LND of horizontal type on $A$, then
$$\ML_{h}(A)=\bigoplus_{m\in \theta\cap\delta_L}\KK\varphi_m\chi^m\,,$$
where $L=L(\partial)$ and $\varphi_m\in A_m$ is such that
$\divi(\varphi_m)+\DD(m)=0$.
\end{theorem}
\begin{proof}
  The assertions follow immediately by virtue of \eqref{ML-eq} and
  Lemmas \ref{ML-fib} and \ref{ML-hor}.
\end{proof}

In the following corollary we give a criterion of triviality of the
homogeneous Makar-Limanov invariant $\ML_h(A)$.
\begin{corollary}\label{ML-triv}
  With the notation as above, $\ML_h(A)=\KK$ if and only if one of the
  following conditions hold.
  \begin{enumerate}[(i)]
  \item $A$ is elliptic, $\rank(M)\geq 2$, and $\deg\DD$ does not
    intersect any extremal ray of $\sigma$.
  \item $A$ admits a homogeneous LND of horizontal type and
    $\theta\cap\delta=\{0\}$.
  \end{enumerate}
  In particular, in both cases $\ML(A)=\KK$.
\end{corollary}

\begin{proof}
  By Lemma \ref{ML-fib}, (i) holds if and only if
  $\ML_{hor}(A)=\KK$. By Theorem \ref{ML-cp1}, (ii) holds if and only
  if there is a homogeneous LND of horizontal type and $\ML_h(A)=\KK$.
\end{proof}

\begin{example}
  It easily seen that $\ML_h(A)=\KK$ for $A$ as in Example
  \ref{ex-hyp-rev}.
\end{example}

\subsection{A non-rational threefold with trivial Makar-Limanov
  invariant}
\label{non-rat}

To exhibit such an example, we let $\sigma$ be a pointed polyhedral
cone in $M_\QQ$, where $\rank(M)=n\geq 2$. We let as before
$A=A[C,\DD]$, where $\DD$ is a proper $\sigma$-polyhedral divisor on a
smooth curve $C$. By \ref{18} $\fract\,A=K_0(M)$ and so $\spec\, A$ is
birational to $C\times \PP^n$ (cf. Corollary 3 in \cite{Ti2}).

By Corollary \ref{ML-triv}, if $A$ is non-elliptic and $\ML(A)=\KK$,
then $A$ admits a homogeneous LND of horizontal type. So
$C\simeq\AF^1$ and $\spec\, A$ is rational. On the other hand, the
curve $C$ does not participate in the assumptions of Corollary
\ref{ML-triv} (i). So if (i) is fulfilled, then $\ML(A)=\KK$ while
$\spec\, A$ is birational to $C\times \PP^n$. This leads to the
following result.

\begin{proposition} \label{non-rat-res} Let $A=A[C,\DD]$, where $\DD$
  is a proper $\sigma$-polyhedral divisor on a smooth projective curve
  $C$ of positive genus. Suppose further that $\deg\DD$ does not
  intersect any extremal ray of $\sigma$. Then $\ML(A)=\KK$ whereas
  $X=\spec\, A$ is non-rational
\end{proposition}

\begin{remark}
  It is evident that $X$ in Proposition \ref{non-rat-res} is in fact
  stably non-rational i.e., $X\times \PP^\ell$ is non-rational for all
  $\ell\geq 0$, cf. \cite[Example~1.22]{Po}.
\end{remark}

In the rest of this section we give a simple geometric example
illustrating this proposition.

Letting $N=\ZZ^2$ and $M=\ZZ^2$ with the canonical bases and duality,
we let $\sigma\subseteq N_\QQ$ be the first quadrant,
$\Delta=(1,1)+\sigma$, and $h=h_{\Delta}$ so that
$h(m_1,m_2)=m_1+m_2$.  Furthermore, we let $A=A[C,\DD]$, where
$C\subseteq \PP^2$ is the elliptic curve with affine equation
$s^2-t^3+t=0$, and $\DD=\Delta\cdot P$ is the proper
$\sigma$-polyhedral divisor on $C$ with $P$ being the point at
infinity of $C$.

Since $C\not\simeq \PP^1$ and $\deg\DD=\Delta$, $A$ satisfies the
assumptions of Corollary \ref{non-rat-res}. Letting $K_0$ be the
function field of $C$, by Theorem \ref{AH-description} we obtain
$$A_{(m_1,m_2)}=H^0(C,\OO_C((m_1+m_2)P))\subseteq K_0\,.$$

The functions $t,s\in K_0$ are regular in the affine part of $C$, and
have poles of order 2 and 3 on $P$, respectively. By the Riemann-Roch
Theorem $\dim H^0(C,\OO(rP))=r$ $\forall r>0$. Hence the functions
$\{t^{i},t^{j}s\mid 2i\leq r\mbox{ and } 2j+3\leq r\}$ form a basis of
$H^0(C,\OO(rP))$ (see \cite{Ha} Chapter IV, Proposition 4.6).

In this setting the first gradded pieces are the $\KK$-modules
$$
\begin{array}{c}
  A_{(0,0)}=A_{(1,0)}=A_{(0,1)}=\KK\,, \\
  A_{(2,0)}=A_{(1,1)}=A_{(0,2)}=\KK\oplus\KK t\,, \\
  A_{(3,0)}=A_{(2,1)}=A_{(1,2)}=A_{(0,3)}=\KK\oplus\KK t\oplus\KK s\,, \\
  A_{(4,0)}=A_{(3,1)}=A_{(2,2)}=A_{(1,3)}=A_{(0,4)}=\KK\oplus\KK t\oplus\KK t^2\oplus\KK s\,.
\end{array}
$$

\begin{remark}
  Let $\EE$ be the locally free sheaf of rank 2
  $\OO_C(P)\oplus\OO_C(P)$. The variety $\spec\, A$ corresponds to the
  contraction of the zero section of the vector bundle associated to
  $\EE$.
\end{remark}

It is easy to see that $A$ admits the following set of generators.
$$
\begin{array}{c}
  u_1=\chi^{(1,0)}, \qquad u_2=\chi^{(0,1)}, \qquad u_3=t\chi^{(2,0)}, \qquad u_4=t\chi^{(1,1)}, \qquad u_5=t\chi^{(0,2)}, \vspace{1.3ex} \\ 
  u_6=s\chi^{(3,0)}, \qquad u_7=s\chi^{(2,1)}, \qquad u_8=s\chi^{(1,2)}, \qquad u_9=s\chi^{(0,3)}\,.
\end{array}
$$
So $A\simeq \KK^{[9]}/I$, where $\KK^{[9]}=\KK[x_1,\ldots,x_9]$, and
$I$ is the ideal of relations of $u_i$ ($i=1\ldots 9$)\footnote{Using
  a software for elimination theory, we were able to find a minimal
  generating set of $I$ consisting of 22 polynomials.}.

Furthermore, $A_m\subseteq \KK[s,t]/(s^2-t^3+t)$ $\forall
m\in\sigma_M^\vee$ since $\DD$ is supported at the point at infinity
$P$. The semigroup $\sigma_M^\vee$ is spanned by $(1,0)$ and $(0,1)$,
so letting $v=\chi^{(1,0)}$ and $w=\chi^{(0,1)}$ we obtain
$$
A=\KK[v,w,tv^2,tvw,tw^2,sv^3,sv^2w,svw^2,sw^3] \subseteq
\KK[s,t,v,w]/(s^2-t^3+t)\,.
$$
Thus $\spec\,A$ is birationally dominated by $C_0\times\AF^2$, where
$C_0=C\setminus\{P\}$.

Since $C\not\simeq \PP^1$, by Lemma \ref{hor-freeness} there is no
homogeneous LND of horizontal type on $A$. There are two extremal rays
$\rho_i\subseteq\sigma$ spanned by the vectors $(1,0)$ and
$(0,1)$. Since $\deg\DD=\Delta$ is contained in the relative interior
of $\sigma$, Corollaries \ref{feq} and \ref{fiex} imply that there are
exactly 2 pairwise non-equivalent homogeneous LNDs $\partial_i$ of
fiber type which correspond to the extremal rays $\rho_i$, $i=1,2$,
respectively.

The codimension 1 face $\tau_1$ dual to $\rho_1$ is spanned by $(0,1)$
and, in the notation of Lemma \ref{fib-core},
$S_{\rho_1}=\{(-1,r)\mid r\geq0\}$. Letting $e_1=(-1,1)$ yields
$D_{e_1}=0$ and so $\Phi_{e_1}=\KK$. We fix $\varphi_1=
1\in\Phi_{e_1}$.  By the same lemma we can chose
$\partial_1=\partial_{\rho_1,e_1,\varphi_1}$ as
$$\partial_1\left(\chi^{(m_1,m_2)}\right)=m_1\cdot\chi^{(m_1-1,m_2+1)},
\quad\mbox{for all}\quad (m_1,m_2)\in \sigma^\vee_M\,.$$

Likewise, the codimension 1 face $\tau_2$ dual to $\rho_2$ is spanned
by $(1,0)$ and, in the notation of Lemma \ref{fib-core},
$S_{\rho_2}=\{(r,-1)\mid r\geq0\}$. Letting $e_2=(1,-1)$ yields
$D_{e_2}=0$ and so $\Phi_{e_2}=\KK$. We fix $\varphi_2=
1\in\Phi_{e_2}$.  By Lemma \ref{fib-core} we can chose
$\partial_2=\partial_{\rho_2,e_2,\varphi_2}$ as
$$\partial_2\left(\chi^{(m_1,m_2)}\right)=m_2\cdot\chi^{(m_1+1,m_2-1)},
\quad\mbox{for all}\quad (m_1,m_2)\in \sigma^\vee_M\,.$$

The kernels of $\partial_1$ and $\partial_2$ are given by
$$\ker\partial_1=\bigoplus_{m\in\tau_1\cap M}A_m\chi^m\quad \mbox{and}\quad \ker\partial_2=\bigoplus_{m\in\tau_2\cap M}A_m\chi^m\,.$$ 
Since $\tau_1\cap\tau_2=\{0\}$ we have
$$\ML(A)=\ker\partial_1\cap\ker\partial_2=A_{(0,0)}=\KK\,.$$
This agrees with Corollary \ref{non-rat-res}.

The LNDs $\partial_i$ are induced, under the isomorphism $A\simeq
\KK^{[9]}/I$, by the following LNDs on $\KK^{[9]}$:
$$
\partial_1=x_2\frac{\partial}{\partial
  x_1}+2x_4\frac{\partial}{\partial x_3}+x_5\frac{\partial}{\partial
  x_4}+3x_7\frac{\partial}{\partial x_6}+2x_8\frac{\partial}{\partial
  x_7}+x_9\frac{\partial}{\partial x_8}\,,$$ and
$$
\partial_2=x_1\frac{\partial}{\partial
  x_2}+x_3\frac{\partial}{\partial x_4}+2x_4\frac{\partial}{\partial
  x_5}+x_6\frac{\partial}{\partial x_7}+2x_7\frac{\partial}{\partial
  x_8}+3x_8\frac{\partial}{\partial x_9}\,,
$$
respectively.

We let below $X=\spec\, A$, and we let $\pi:X\dashrightarrow C$ be the
rational quotient for the $\TT$-action on $X$. The comorphism of $\pi$
is given by the inclusion $\pi^*:K_0\hookrightarrow
\fract\,A=K_0(u_1,u_2)$.

The orbit closure $\Theta=\overline{\pi^{-1}(0,0)}$ over $(0,0)\in C$
is general and it is isomorphic to $\AF^2=\spec\,\KK[x_1,x_2]$. The
restrictions to $\Theta$ of the $\KK_+$-actions $\phi_i$ corresponding
to $\partial_i$, $i=1,2$, respectively are given by
$$\phi_1|_\Theta:(t,(x_1,x_2))\mapsto(x_1+tx_2,x_2)\quad \mbox{and} \quad
\phi_2|_\Theta:(t,(x_1,x_2))\mapsto(x_1,x_2+tx_1)\,.$$

Furthermore, there is a unique singular point $\bar{0}\in X$
corresponding to the fixed point of the $\TT$-action on $X$. The point
$\bar{0}$ is given by the augmentation
ideal $$A_+=\bigoplus_{\sigma_M^\vee\setminus\{0\}}A_m\chi^m\,,$$

On the other hand, let $A=A[C,\DD]$, where $\DD$ is a proper
$\sigma$-polyhedral divisor on a smooth projective curve $C$. By
Theorem 2.5 in \cite{KaRu}, if $\spec\, A$ is smooth, then $\spec\,
A\simeq\AF^{n+1}$ (see also Proposition 3.1 in \cite{Su}). In
particular, $\spec\, A$ is rational.

\begin{remark}
  \begin{enumerate}[$(i)$]
  \item In \cite{Li1} generalizing the methods of this section we
    obtain a birational characterization of normal affine varieties
    with trivial ML-invariant.

  \item In \cite{Li2} we studied singularities of $\TT$-varieties. In
    particular, we showed that the singularities of the $X=\spec\,
    A[C,\DD]$ are not Cohen-Macaulay. On the other hand, in the recent
    preprint \cite{Po} a new family of examples of non-rational affine
    varieties with trivial ML-invariant is given. This time, these
    varieties are smooth.
  \end{enumerate}
\end{remark}

\end{document}